\documentclass[a4paper,11pt]{amsart}

\usepackage[all]{xy}
\usepackage{amssymb,amsmath,eufrak}
\author[Florent Benaych-Georges]{Florent Benaych-Georges}

\address{Florent Benaych-Georges, LPMA, UPMC Univ Paris 6, 4, Place Jussieu, 75005 Paris, France.}
\title[Cycles of free words in random permutations]{Cycles of free words in several independent random permutations with restricted cycle lengths}
\date{\today}
\newcommand{\Sy}{\mathfrak{S}}
\newcommand{\Po}{\operatorname{P}}

\newcommand{\Poisson}{\operatorname{Poisson}}
\newcommand{\Poiss}{\operatorname{Poisson}}
\newcommand{\Beg}{\operatorname{Beg}}
\newcommand{\End}{\operatorname{End}}

\newcommand{\Part}{\operatorname{Part}}

\newcommand{\Adm}{\operatorname{Adm}}

\newcommand{\ds}{\displaystyle}

\newcommand{\ssi}{if and only if }

\newcommand{\E}{\mathbb{E}}

\newcommand{\pro}{probability }

\newcommand{\f}{\frac}
\newcommand{\ff}{\frac{1}}
\newcommand{\lf}{\left}
\newcommand{\ri}{\right}

\newcommand{\st}{such that }
\newcommand{\la}{\lambda}

\newcommand{\tii}{\scriptstyle\times\ds\!}

\newcommand{\vfi}{\varphi}
\newcommand{\ste}{\, ;\, }
\newcommand{\mc}{\mathcal }

\newcommand{\eps}{\varepsilon}

\newcommand{\bck}{\backslash}
\newcommand{\si}{\sigma}

\newtheorem{Th}{Theorem}[section]
\newtheorem{propo}[Th]{Proposition} 
 
\newtheorem{lem}[Th]{Lemma}

\newtheorem{rmq}[Th]{Remark}
\newtheorem{cor}[Th]{Corollary}

\newenvironment{pr}{\noindent {\bf Proof. }}{\hfill $\square$}

\newenvironment{prlem}{\noindent {\bf Proof of the lemma. }}{ \hfill  $\square$}
\long\def\symbolfootnote[#1]#2{\begingroup
\def\thefootnote{\fnsymbol{footnote}}\footnote[#1]{#2}\endgroup}

\begin{document}
\maketitle

 

\begin{abstract}
In this text, 
we consider random permutations which can be written as free words in several independent random permutations: firstly, we fix  a non trivial word $w$  in letters 
$g_1,g_1^{-1},\ldots, g_k,g_k^{-1}$, secondly, for all $n$, we introduce a $k$-tuple $s_1(n),\ldots, s_k(n)$ of independent random permutations of $\{1,\ldots, n\}$, and the random permutation $\sigma_n$ we are going to consider is the one obtained by replacing each letter $g_i$ in  $w$ by $s_i(n)$. For example, for $w=g_1g_2g_3g_2^{-1}$,  $\sigma_n=s_1(n)\circ s_2(n)\circ s_3(n)\circ s_2(n)^{-1}$.  Moreover, we  restrict  the set of possible lengths of the cycles of the $s_i(n)$'s: we fix sets $A_1,\ldots, A_k$ of positive integers and suppose that for all $n$, for all $i$, $s_i(n)$ is uniformly distributed on the set of permutations of $\{1,\ldots, n\}$  which have all their cycle lengths in  $A_i$. 
For all positive integer $l$, we are   going to give asymptotics, as $n$ goes to infinity, on the number $N_l(\sigma_n)$  of cycles of length $l$ of $\sigma_n$. We shall also  consider the joint distribution of the random vectors $(N_1(\sigma_n),\ldots, N_l(\sigma_n))$. 
We first prove that the   of $w$ in a certain quotient of the free group with generators $g_1,\ldots, g_k$ determines the rate of growth of the random variables $N_l(\sigma_n)$ as $n$ goes to infinity. We also prove that in many cases, the distribution of $N_l(\sigma_n)$ converges to a Poisson law with parameter $1/l$ and that the  random variables $N_1(\sigma_n),N_2(\sigma_n), \ldots$ are asymptotically independent.  We notice the   surprising fact that from this point of view, many things happen as if $\sigma_n$ were uniformly distributed  on the $n$-th symmetric group. 
\end{abstract}

\section{Introduction}\subsection{General introduction} 
It is well known \cite{bat} that if $\sigma_n$ is a uniformly distributed random permutation of $\{1,\ldots, n\}$ and for all $l\geq 1$, one denotes the number of cycles of $\sigma_n$ with length $l$  by $N_l(\sigma_n)$, then for all $k\geq 1$,  the joint distribution of $(N_1(\sigma_n), \ldots,N_k(\sigma_n) )$ converges in distribution, as $n$ goes to infinity, to $\Poiss (1/1)\otimes \cdots\otimes \Poiss (1/k)$. The law of ``large cycles" can be recovered using the  Poisson-Dirichlet distribution, but in this paper, we are only going to deal with ``small cycles", i.e. cycles with a size which does not vary with $n$. 

In \cite{nica94}, Nica was the first to consider the case where $\sigma_n$ is a free word in several independent random permutations. For a fixed $k\geq 1$ and a fixed free word $w$ in the letters $g_1,g_1^{-1},\ldots, g_k,g_k^{-1}$, he introduced, for all $n$, a random permutation $\si_n$  of $\{1, \ldots, n\}$ defined in the following way: $\sigma_n$ is the random permutation obtained by replacing each letter $g_i$ in  $w$ by $s_i(n)$, where  $s_1(n)$, \ldots, $s_k(n)$ are independent random permutations uniformly distributed on the $n$-th symmetric group. He proved that under the hypothesis that  $w$ is not a power of another word,  for all positive integer $l$,  $N_l(\sigma_n)$ converges in distribution, as $n$ goes to infinity, to a Poisson distribution with mean $1/l$. His result was reproved recently by Linial and Puder in \cite{lp08}, where it appeared that free words in random permutations are relevant to analyze spectrums of $n$-lifts of graphs.  None of   these papers considered the joint distribution of the $N_l(\sigma_n)$'s.

Independently, in \cite{yakymiv} and in \cite{fbg.1},  the authors  considered random permutations whose distribution is not uniform on the whole symmetric group, but on the set of permutations of $\{1,\ldots, n\}$ which have all their cycle lengths in a fixed set $A$ of positive integers. This kind of random permutation has been the subject of study by a number of authors for more than thirty 
years (for a review of the literature in this area, see \cite{yakymiv1}).  It was proved  \cite{fbg.1} that if $A$ is finite with greatest element $d$, then as $n$ goes to infinity, such a random permutation tends to get not  far away from having   $d$: the cardinality of the subset of $\{1,\ldots, n\}$ covered by the supports of cycles with length $d$ in such a random permutation is asymptotic to $n$. It was also proved \cite{yakymiv,fbg.1} that if $A$ is infinite, then the order of   such a random permutation in the symmetric group  goes to infinity as $n$ does.

In the present paper, we shall somehow mix both of the previously presented  problems: we shall fix $k\geq 1$, a free word $w$ in the letters $g_1,g_1^{-1},\ldots, g_k, g_k^{-1}$, and a list $A_1,\ldots, A_k$ of sets of  positive integers satisfying a technical assumption: \begin{equation}\label{tech.as.14.12.07} \forall i,\quad\quad  \textrm{$A_i$ is finite \quad or \quad satisfies $\sum_{j\geq 1, j\notin A_i} \ff{j}<\infty$}.\end{equation} Then for all $n$ such that it is possible, we shall introduce a $k$-tuple   $s_1(n)$, \ldots, $s_k(n)$ of independent random permutations such that for all $i$, $s_i(n)$ is uniformly distributed on the set of permutations of $\{1,\ldots, n\}$ having all of their cycle lengths in $A_i$ and define $\sigma_n$ to be the random permutation obtained by replacing each letter $g_i$ in $w$ by $s_i(n)$.  

The starting point of this study was the  relation between random permutations  and asymptotic first and second order freeness. The {\it asymptotic freeness} of   random matrices is a notion due  to Voiculescu based upon  the {\it non-commutative distribution} of the matrices (which is, roughly, defined by the (random) algebra generated by the matrices, up to a global conjugation by a unitary matrix). One mights find several introductions to this theory in  \cite{vdn91,hiai,ns06,agz09}. A complementary theory, the theory of  {\it second order freeness}, has been developped these last five years about Gaussian fluctuations of the non-commutative distribution of asymptotically free random matrices around its limit (see 
  \cite{mingo-nica04, mingo-speicher06, mingo-piotr-speicher07, mingo-piotr-collins-speicher07}).    In  \cite{nica93,neagu05} the authors proved that   random matrices associated with uniform random permutations (with, possibly, restricted cycle length) are asymptotically free, i.e. that, as far as their joint non-commutative distribution is concerned, they behave like large Haar-distributed unitary matrices. Concretely, it means that the number of fixed points  any non-empty word in such random permutations is $o(n)$ (the fact that this is true for {\it any} non-empty word is of importance here). The following step, in the study of this issue, was to consider  the asymptotic distribution of the number of fixed points of free words in random permutations (with, possibly, restricted cycle length). In this paper, we extend the investigation to the number of cycles of any fixed length, and prove that for several models (depending on the way we restrict the cycle length), the asymptotic distribution is  Poisson.

A more detailed review of our results is given below. 
 \subsection{Link with a quotient of the free group generated by $g_1,\ldots, g_k$}\label{bobo.06.07}Since, as explained above, for $i$ \st $A_i$ is finite, $s_i(n)$ is not far away from having order $d_i:=\sup A_i$, 
 it is natural to expect that for large values of $n$,  the distribution of $\sigma_n$ will depend firstly on the word obtained from $w$ by removing all sequences of the type $g_i^{\pm d_i}$ for $i$ \st $A_i$ is finite. This is what has led us to introduce the group $F_k/[g_1^{d_1},\ldots, g_k^{d_k}]$ generated by free elements $g_1,\ldots, g_k$ and quotiented by the relations $g_i^{d_i}=1$ for $i\in [k]$, with $d_i=\sup A_i$ (when $d_i=+\infty$, the relation $g_i^{d_i}=1$ is not included). Our more general results, respectively Theorems  \ref{15.6.07.1} and \ref{15.11.06.989}, are the followings: \begin{itemize} \item[-] If the element of $F_k/[g_1^{d_1},\ldots, g_k^{d_k}]$ represented by $w$ has finite order $d\geq 1$, then, as $n$ tends to infinity,  $N_d(\sigma_n)/n$ converges to $1/d$  and for all $l\neq d$,  $N_{l}(\sigma_n)/n$ converges to $0$. This means that  $\sigma_n$ is not far away from having order $d$: the cardinality of the subset of $\{1,\ldots, n\}$ covered by the supports of the cycles with length $d$ in such  a random permutation  is asymptotic to $n$.
\item[-] If  the element of $F_k/[g_1^{d_1},\ldots, g_k^{d_k}]$ represented by $w$ has infinite order, then  two cases can occur:
  \begin{itemize}
\item[(i)]  
    There is $i\in [k]$ \st $A_i$ is infinite and $\alpha$ a nonzero integer \st $w$ represents, up to conjugation,  the same element, in $F_k/[g_1^{d_1},\ldots, g_k^{d_k}]$, as the word    $g_i^\alpha$: then for all $l\geq 1$ \st $l|\alpha|\in A_i$,  $$\ds\liminf_{n\to\infty} (\E(N_l(\sigma_n)))\;\geq\; 1/l.$$
    \item[(ii)]  The element of  $F_k/[g_1^{d_1},\ldots, g_k^{d_k}]$ represented by $w$ is not conjugate to an element represented by a word of the type $g_i^\alpha$, with $i\in [k]$, $\alpha$ an integer: then for all $l\geq 1$, 
     $$\ds\liminf_{n\to\infty}( \E(N_l(\sigma_n)))\;\geq \;1/l.$$
\end{itemize}
The result  (ii)  means that  the cycles  of the letters of $w$ are going to mix sufficiently well  to give rise to cycles of most lengths, at least as much as in a uniform random permutation, even when the letters of $w$ have very restricted cycle lengths.  \end{itemize} 

\subsection{Weak limit of the distributions of the $N_l(\sigma_n)$'s}\label{23h43.15.11.06} Then, we are going to prove more precise results:  under certain hypotheses on $w$, for all $l\geq 1$, the joint distribution of the random vector $$(N_1(\sigma_n),\ldots, N_l(\sigma_n))$$ converges weakly, as $n$ goes to infinity, to $$\Poisson(1/1)\otimes\Poisson (1/2)\otimes \cdots \otimes \Poisson (1/l),$$ just as if the distribution of $\sigma_n$ would have been uniform on $\Sy_n$. When all $A_i$'s are infinite (Theorem \ref{2046.15.11.06.1}), our hypotheses  are that  $w$ is neither a single letter nor a power of another word. 
If some of the $A_i$'s are finite, we prove  this result for $w=g_1\cdots g_k$ (Theorem \ref{downpressure.peter.teuchi}), with the exception of the case $k=2$ and $A_1\cup A_2\subset \{1,2\}$ (i.e. of the product of two random involutions), where  we prove that  the $N_l(\sigma_n)$'s are still asymptotically independent, but that their asymptotic distributions are slightly different (Theorem \ref{mardi.28.11.06.2}).

These results extend on the one hand a result of Nica (reproved by Linial and Puder in \cite{lp08}), who proved in \cite{nica94} that in the case where all $A_i$'s are equal to the set of all positive integer and where $w$ 
is not a power of another word, for all $l\geq 1$, the distribution of $N_l(\sigma_n)$ converges weakly to $\Poiss (1/l)$ (without considering the joint distribution), and on the other hand a result of  Neagu, who proved in \cite{neagu05} that  the number of fixed points of $\sigma_n$ is $o(n)$. Here, we shall mention that much of the methods we use in this paper 	are inspired by the ones introduced in \cite{nica94,neagu05}, and that we use many of their results.

\subsection{Comments on these results and open questions}

a) In this paper, only small cycles (i.e. cycles which size does not depend on $n$)  are considered by our investigations. It would be interesting to know if it is possible to prove a limit theorem for the whole random partition of $\{1,\ldots,n\}$ defined by the cycles of $\sigma_n$. It is natural to expect a Poisson-Dirichlet distribution \cite{pitman,bat}.

b)   Another   interesting question arising from the similarity between $\sigma_n$ and a uniform random permutation on $\Sy_n$ is the following: do we have a characterization of the words $w$ \st  for $n$ large enough, $\sigma_n$  is uniformly distributed ?

\subsection{Notation} In this text, for $n$ an integer, we shall denote by $[n]$ the set $\{1,\ldots, n\}$ and   by $\Sy_n$ the group of permutations of $[n]$. For $A$ a set of positive integers, $\Sy_n(A)$ denotes  the set of permutations of $[n]$ whose cycles have length in $A$.   
 For $\sigma\in \Sy_n$ and $l\geq 1$, we shall denote by $N_l(\sigma)$ the number of cycles of length $l$ in the decomposition of $\sigma$ as a product of cycles with disjoint supports. For $\la >0$, $\Poiss(\la)$ will denote the Poisson distribution with parameter $\la$.

\section{Combinatorial preliminaries to the study of  words in random permutations}
\subsection{Words and groups generated by relations}\label{sectionfreewords.12.10.06}
\subsubsection{Words}\label{Iwashappyinthehazeofadrunkenhour.words}
Let, for $k\geq 1$,  $\mathbb{M}_k$ be the set of {\it words} in the letters $g_1,g_1^{-1}, \ldots, g_k,g_k^{-1}$, i.e. the set of sequences $g_{i_1}^{\alpha_1}\cdots g_{i_{n}}^{\alpha_{n}}$, with $n\geq 0$,  $i_1,\ldots, i_{n}\in [k]$, $\alpha_1,\ldots, \alpha_{n}=\pm 1$. A word $w\in \mathbb{M}_k$ is said to be {\it reduced} if in its writing, no letter is followed by its inverse. It is said to be {\it cyclically reduced} if moreover, the first and the last letters are not the inverses of each other.  


 A cyclically reduced word is said to be {\it primitive} if it is not the concatenation of $d\geq 2$ times the same word. 
 

For $w=g_{i_1}^{\alpha_1}\cdots g_{i_{|w|}}^{\alpha_{|w|}}\in \mathbb{M}_k $ and $s=(s_1,\ldots, s_k)$ a $k$-tuple of elements of a group, $w(s)$ denotes  $s_{i_1}^{\alpha_1}\cdots s_{i_{|w|}}^{\alpha_{|w|}}$.

\subsubsection{The quotient of the free group with $k$ generators by the relations $g_1^{d_1}=1,\ldots, g_k^{d_k}=1$}\label{sectionfreewordsquotiented.12.10.06} Let $F_k$ be the free group generated by $g_1, \ldots, g_k$. This is the set of reduced words of $\mathbb{M}_k$ endowed with the operation of concatenation-reduction via the relations $g_ig_i^{-1}=1,  g_i^{-1}g_i=1, i\in [k]$. For $w, w_0\in \mathbb{M}_k$ with $w_0$ reduced, $w$ is said to {\it represent} or to be a {\it writing of} the element $w_0$ of $F_k$  if one can reduce $w$ to $w_0$ via the previous relations.  
For example,  $g_1g_3g_3^{-1}g_1$ is a writing of $g_1^2$.

Consider  $d_1,\ldots, d_k\in\{2,3,4,\ldots\}\cup\{\infty\}$. Let $F_k/[g_1^{d_1},\ldots, g_k^{d_k}]$ be the group $F_k$ quotiented by its normal subgroup generated by the set $\{g_i^{d_i}\ste i\in [k], d_i<\infty\}$.  A word $w\in \mathbb{M}_k$ is said to {\it represent} or to be a {\it writing of} an element  $C$ of $F_k/[g_1^{d_1},\ldots, g_k^{d_k}]$  if it is a writing of an element of $C$ (seen as a subset of $F_k$).  

Theorem 1.4 of Section  1.4 of \cite{MKS} states the following facts. 
\begin{Th}\label{atalanta-bergame.11.06}(a) Any element of $F_k/[g_1^{d_1},\ldots, g_k^{d_k}]$ has a writing of the type $g_{i_1}^{\alpha_1}\cdots g_{i_{n}}^{\alpha_{n}}$
 with $n\geq 1$,  $i_1\neq i_2\neq \cdots \neq i_n\in [k]$, $ \alpha_1,\ldots, \alpha_n$ integers \st  $0<|\alpha_1|<d_{i_1}$,\ldots,  $0<|\alpha_n|<d_{i_n}$, and this writing is unique up to replacements of the type  $$g_{i_j}^{\alpha_j}\to \begin{cases}g_{i_j}^{d_{i_j}+\alpha_j} & \textrm{if $\alpha_j<0$ and $d_{i_j}<\infty$,}\\
 g_{i_j}^{-d_{i_j}+\alpha_j} & \textrm{if $\alpha_j>0$ and $d_{i_j}<\infty$,}
 \end{cases}$$
 with $j\in [n]$. 
 
 (b) In any conjugacy class of  the group $F_k/[g_1^{d_1},\ldots, g_k^{d_k}]$, there is an element represented by a word  of the previous type \st moreover, $i_n\neq i_1$, and such a word is  unique up to replacements of the previous  type and to transformations of the type 
$g_{i_1}^{\alpha_1}\cdots g_{i_{n}}^{\alpha_{n}}\to g_{i_{n}}^{\alpha_{n}}g_{i_1}^{\alpha_1}\cdots g_{i_{n-1}}^{\alpha_{n-1}}$. 
\end{Th}


Let us define the {\it $(d_1 \ldots, d_k)$-cyclically reduced words} to be the  words of the type $g_{i_1}^{\alpha_1}\cdots g_{i_{n}}^{\alpha_{n}}$
 with $n\geq 0$, $i_1\neq i_2\neq \cdots \neq i_n\neq i_1\in [k]$, $ \alpha_1,\ldots, \alpha_n$ integers \st  $0<|\alpha_1|<d_{i_1}$,\ldots,  $0<|\alpha_n|<d_{i_n}$. 
 
 
 Consider a cyclically reduced word $w\in \mathbb{M}_k$.  We shall call a {\it partial  $(d_1,\ldots, d_k)$-cyclic reduction of $w$} a word which can be obtained from $w$ in a finite number of steps of the following types:
\begin{itemize}
\item[(a)] $uv\to vu$, with  $u,v\in \mathbb{M}_k$,
\item[(b)] $ug_i^{\alpha d_i}u^{-1}v\to v$, with $i\in [k]$, $\alpha$ a nonzero integer, $u$ a nonempty word and $v\in \mathbb{M}_k$ a cyclically reduced word,
\item[(c)] $g_i^{\alpha d_i}v\to v$, with $i\in [k]$, $\alpha$ a nonzero integer, and  $v$ a word which can be written   $v=g_i^{\beta}v'$, with $\beta$ an integer, $|\beta|<d_i$ and $v'$  a word whose first and last letters do not belong to $\{g_i,g_i^{-1}\}$.
\end{itemize} 
A partial $(d_1,\ldots, d_k)$-cyclic reduction of $w$ which is $(d_1 \ldots, d_k)$-cyclically reduced is  called a {\it $(d_1 \ldots, d_k)$-cyclic reduction} of $w$. In general, there is more than one $(d_1 \ldots, d_k)$-cyclic reduction of $w$, but if there is only one, we call it {\it the} $(d_1 \ldots, d_k)$-cyclic reduction of $w$.

 For example,  for $d_1=4, d_2=5$,  the following words are partial $(d_1,\ldots, d_k)$-cyclic reductions of $w=g_2^2g_1g_2^{6}g_3g_1^{-4}g_3^{-1}g_1^{-1}g_2^3$:  $$g_2^5g_1g_2^{6}g_3g_1^{-4}g_3^{-1}g_1^{-1}, \;\;\,\;g_2^{6}g_3g_1^{-4}g_3^{-1},\, \;\; \;g_2g_3g_1^{-4}g_3^{-1},\, \;  \;\;g_3g_1^{-4}g_3^{-1}g_2,\,\; \;\;
  g_2.$$ The last one is a  $(d_1 \ldots, d_k)$-cyclic reduction of $w$. 

\subsection{Admissible graphs and partitions. Colored graphs associated to words and permutations}\label{combi.prelim.misssyou.cecilegoestoustomorrow}
\subsubsection{Basic graph theoretic  definitions}\label{Lolene-Yola}
 In this text, we shall consider {\it oriented edge-colored  graphs with color set $[k]$}. These are families $G=(V; E_1,\ldots, E_k)$, where $V$ is a finite set  (its elements are called the {\it vertices} of $G$) and for all $r\in [k]$, $E_r$ is a subset of $V^2$ (the set of {\it edges with color $r$} of $G$).  
 For $e=(u,v)$ edge of $G$, $u$ (resp. $v$) is called  the  {\it   beginning vertex} of $e$ (resp. the   {\it ending vertex} of $e$) and is denoted by $\Beg(e)$ (resp. $\End(e)$). $e$ is often denoted by $u\to v$. Note that  the $E_r$'s are not supposed to be pairwise disjoint.   
  Throughout this paper, the color set of the edges will always be $[k]$, so it will often be implicit. 
 
For $r\in [k]$, $G=(V; E_1,\ldots, E_k)$ is  said to be {\it monochromatic} with {\it color} $r$ if for all $i\neq r$, $E_i=\emptyset$. If moreover, $$
V= \{v_1,\ldots,v_l\}\textrm{ and }E_r=\{v_1\to v_2\to\cdots\to v_l\}\textrm{ (resp. $E_r=\{v_1\to v_2\to\cdots\to v_l\to v_1\}$)},$$  with $l\geq 1$ and $v_1,\ldots, v_l$ pairwise distinct, $G$ is said to be a   {\it monochromatic directed path} (resp. {\it monochromatic directed cycle}) with    {\it length}  $l-1$ (resp. $l$).

A {\it subgraph} of a graph $(V; E_1, \ldots, E_k)$ is a graph of the type $(W; F_1, \ldots, F_k)$, with $W\subset V, F_1\subset E_1,\ldots, F_k\subset E_k.$

A {\it monochromatic directed cycle} (resp. a {\it monochromatic directed path}) of an oriented edge-colored graph is a subgraph which is a monochromatic directed cycle (resp. a monochromatic directed path). 
 
 A graph $(V; E_1, \ldots, E_k)$ is said to be the {\it disjoint union} of the graphs  $(V'; E'_1, \ldots, E'_k)$ and $(V''; E''_1, \ldots, E''_k)$ if $V=V'\cup V''$, $E_1=E_1'\cup E_1''$, \ldots, $E_k=E_k'\cup E_k''$ and all of these unions are disjoint.







An {\it isomorphism} between two oriented edge-colored  graphs $(V; E_1, \ldots, E_k)$, $(V'; E'_1, \ldots, E'_k)$ is a bijection $\vfi : V\to V'$ \st for all $u,v\in V$, for all $r\in [k]$, one has $$u\to v \in E_r\Longleftrightarrow \vfi(u)\to \vfi(v)\in E'_r.$$

\subsubsection{Admissible graphs and partitions}

Recall that a {\it partition} $\Delta$ of a set $X$ is a set of pairwise disjoint, nonempty subsets of $X$  (called the {\it classes} of $\Delta$) whose union is $X$. In this case, for $x,y\in X$,  ``$x=y\mod \Delta$" and ``$\Delta$ links $x$ and $y$" both mean that $x,y$ are in the same class of $\Delta$. Since $\Delta$ is a set, $|\Delta|$ denotes  its cardinality. 

For any function $\gamma$ defined on a set $X$, we shall denote by $\Part(\gamma)$ the partition of $X$ by the level sets of $\gamma$.

Let $G$ be an oriented edge-colored graph with color set $[k]$  and with vertices set $V$.

 $G$ is said to be {\it admissible} if two different edges with the same color cannot have the same beginning or the same ending vertex.

We are going to use the notion of {\it quotient graph}, that we define now. Let us define, for $\Delta $ a partition of $V$, $G/\Delta$ to be  the oriented edge-colored graph    whose vertices are the classes of $\Delta$ and such that for all $C,C'$ classes of  $\Delta$, for all $r\in [k]$, there is an  edge with color $r$ from $C$ to $C'$ in $G/\Delta$ when there is an edge with color $r$ in $G$ from a vertex of $C$ to a vertex of $C'$. 

\begin{rmq}\label{9.05.07.1}Note that clearly, if $\Delta_1$ is a partition of $V$ and $\Delta_2$  is a partition of $\Delta_1$, then the oriented edge-colored graph $(G/\Delta_1)/\Delta_2$ is isomorphic to $G/\Gamma$, where $\Gamma$ is the partition of $V$ \st for all $x,y\in V$, $x=y\mod \Gamma$ if and only if the classes of $x$ and $y$ in $\Delta_1$ are in the same class of $\Delta_2$.
\end{rmq}

A partition $\Delta$ of $V$  is said to be an {\it admissible partition of $G$} when  the graph $G/\Delta$ is admissible, i.e.  when for all pair $(e,f)$ of edges of $G$ with the same color, $$\Beg(e)=\Beg(f)\mod \Delta\iff \End(e)=\End(f)\mod \Delta.$$  

In the following proposition, whose proof is obvious, we define the operator $\Adm$ on the set of  oriented edge-colored graphs.

\begin{propo}Let $G$ be an oriented edge-colored graph. Let $\Delta$ be the partition of the vertex set of $G$ which links two vertices $r,s$ \ssi there are $n\geq 0$,  $t_0=r, t_1, \ldots, t_{n}=s$ some vertices of $G$,  $i_1,\ldots, i_n\in [k]$, and  $\eps_1,\ldots, \eps_n\in \{\pm 1\}$ such that:
\\
 - the reduction of the word $g_{i_1}^{\eps_1}\cdots g_{i_n}^{\eps_n}$ is the empty word,
 \\
 - for all $l\in [n]$, the $i_l$-colored edge $t_{l-1}\to t_l$ (resp. $t_{l-1}\leftarrow t_l$) belongs to $G$ if $\eps_l=1$ (resp. if $\eps_l=-1$).

Then $\Delta$ is the minimal partition $P$ (with respect to the refinement order)  \st $G/P$ is admissible.
\end{propo}
With the notations of the previous proposition, the partition $\Delta$ will be called the {\it minimal admissible partition} of $G$ and $G/\Delta$ will be denoted by $\operatorname{Adm}(G)$.

\begin{rmq}\label{ext<->adm}Let $H$ be an oriented edge-colored graph with vertex set $W$ and let $G$ be a subgraph of $H$ with vertex set $V$. Then  $\Adm(H)$ can be obtained (up to an isomorphism) from $\Adm(G)$, by adding to $\Adm(G)$  the vertices and edges of $H$ which are not in $G$, and then quotienting the obtained graph by its minimal admissible partition. 
More specifically, if, for each vertex $v$ of $G$,  one denotes by $\overline{v}$ the class of $v$ in the minimal admissible partition of $G$, then $\Adm(H)$ can be identified, via an isomorphism, with $\Adm(H')$, where $H'$ is the oriented edge-colored graph whose vertex set is the union of the vertex set of $\Adm(G)$ with the set $W\bck V$ 
and whose set of edges of color $i$, for each $i\in [k]$, is the union of the set of  edges of color $i$ of $\Adm(G)$ with the set $$\{w\to w'\ste w,w'\in W\bck V, w\to w'\textrm{ $i$-colored edge of $H$}\}$$ $$\cup\{\overline{v}\to w\ste v\in V, w\in W\bck V,  v\to w\textrm{ $i$-colored edge of $H$}\}$$ $$\cup\{w\to \overline{v}\ste v\in V, w\in W\bck V, w\to v\textrm{ $i$-colored edge of $H$}\}.$$
\end{rmq}

\begin{lem}\label{Adm(ext)=Adm(quotient)}Consider an oriented edge-colored graph $G=(V;E_1,\ldots, E_k)$, $r,s,t\in V$, $i\in[k]$ \st $s\to t\in E_i$. Then for $H$ the  oriented edge-colored graph obtained from $G$ by adding, if it is not already in $G$, the $i$-colored edge $r\to t$ (resp. $s\to r$), $\Adm(H)$ is isomorphic to $\Adm(G/\{r=s\})$ (resp. $\Adm(G/\{r=t\})$), where for $u,v\in V$, $\{u=v\}$ denotes the partition of $V$ whose classes are all singletons, except one: $\{u,v\}$.
\end{lem}

\begin{pr}The proof is immediate with Remark \ref{9.05.07.1} and the following observation:  for $\Gamma$ an admissible partition of $G$, if, for any $v\in V$,  one denotes  by $\overline{v}$  the class of $v$ in $\Gamma$, then  $$\overline{r}=\overline{s}\textrm{ (resp.  $\overline{r}=\overline{t}$)}\iff \overline{r}\to\overline{t}\textrm{ (resp. $\overline{s}\to\overline{r}$) is an $i$-colored edge of $G/\Gamma$}. $$
 \end{pr}

\subsubsection{The graph $G(\sigma,w)$}\label{def.de.H*w} 
Fix $k,p\geq 1$, $\sigma\in \Sy_p$ and $w=g_{i_1}^{\alpha_1}\cdots g_{i_{|w|}}^{\alpha_{|w|}}\in \mathbb{M}_k$ a nonempty word. In \cite{nica94}, Nica defined $G(\sigma,w)$ (denoted by $\mc{H}_{\sigma^{-1}}\bigstar w$ in his paper) to be the directed, edge-colored 
  graph with vertex set $V:=[p]\times [|w|]$ and whose edges are the followings:
for all $(m, l)\in [p]\times [|w|]$, \begin{eqnarray*}\alpha_l=1\;&\Rightarrow&\; (m,l)\to \begin{cases}(m,l+1)&\textrm{if $l\neq |w|$,}\\ (\sigma^{-1}(m),1)&\textrm{if $l=|w|$,}\end{cases}\textrm{ is an $i_l$-colored edge of $G(\sigma,w)$},\\ 
\alpha_l=-1\;&\Rightarrow&\; (m,l)\leftarrow \begin{cases}(m,l+1)&\textrm{if $l\neq |w|$,}\\ (\sigma^{-1}(m),1)&\textrm{if $l=|w|$,}\end{cases}\textrm{ is an $i_l$-colored edge of $G(\sigma,w)$.}
\end{eqnarray*}


In the case where $p=1$ and $\sigma=Id$, we shall denote $G(\sigma, w)$ by $G(w)$ and identify its vertex set  with $[|w|]$ by $(1, l)\simeq l$  for all $l\in [|w|]$. We will also use the convention that if $w$ is the empty word, then $G(w)$ is the graph having $1$ for only vertex and no edge.

For example, for $w=g_1g_2g_3g_4g_2^{-1}g_1g_2^{-1}g_5$, $G(w)$ is the graph 
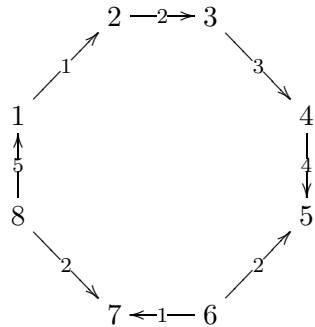
\begin{displaymath} 
\xymatrix{ &2\ar[r]|2&3\ar[dr]|3\\
1\ar[ur]|1&&&4\ar[d]|4\\ 
8\ar[dr]|2\ar[u]|5&&&5\\
& 7&6\ar[l]|1\ar[ur]|2
 } 
\end{displaymath} 
where the colors  of the edges appear on them. We shall give examples of graphs of the type $G(\sigma,w)$ after the following lemma, which expresses $G(\sigma, w)$ as a disjoint union of graphs isomorphic to $G(w^\alpha)$, with $\alpha\geq 1$.

\begin{lem}\label{decomp.de.base} $G(\sigma, w)$ is the disjoint union of $N_1(\sigma) $ graphs isomorphic to $G(w^1)$, $N_2(\sigma)$ graphs isomorphic to $G(w^2)$, \ldots,  $N_p(\sigma)$ graphs isomorphic to $G(w^p)$.
\end{lem}

\begin{pr} Let us denote $G(\sigma, w)=(V; E_1,\ldots, E_k)$. If $I,J$ are disjoint subsets of $[p]$ stabilized by $\sigma$, then in $G(\sigma, w)$, there is no edge between elements of $I\tii [|w|]$ and $J\tii [|w|]$.  Hence $G(\sigma, w)$ is the disjoint union of the graphs $(V(c); E_1(c),\ldots, E_k(c))$, where $c$ varies over the set of cycles of $\sigma$, and where for any such cycle $c$, with support $C\subset [p]$, $V(c)=C\tii [|w|]$ and for any $i\in [k]$, $E_i(c)=E_i\cap (V(c)^2)$.

Hence it suffices to prove that for any $d\in [p]$, for any cycle $c=(m_1\,m_2\cdots m_d)$ of $\sigma$ of length $d$,   there is an isomorphism between $(V(c); E_1(c),\ldots, E_k(c))$ and $G(w^d)$.  The function which maps 
$(m_i, l)\in \{m_1,m_2,\ldots, m_d\}\tii [|w|]$ to $(d-i)|w|+l\in [d|w|]$
is such an isomorphism. 
 \end{pr}

For example, when $w=g_1g_2g_1^{-1}g_2^{-1}$, $p=3$
and $\sigma$ is the cycle $(123)$, $G(\sigma,w)$ is the graph 
$$\xymatrix{ 
(1,1)\ar[r]|1\ar[d]|2& (1,2)\ar[r]|2&(1,3)&(1,4)\ar[l]|1&(3,1)\ar[l]|2\ar[r]|1&(3,2)\ar[d]|2\\
(2,4)\ar[r]|1&(2,3)&(2,2)\ar[l]|2&(2,1)\ar[l]|1\ar[r]|2&(3,4)\ar[r]|1&(3,3)
 } 
$$
where the colors  of the edges appear on them. 

When $w$ is still $g_1g_2g_1^{-1}g_2^{-1}$, but $p=5$
and $\sigma$ is the product of disjoint cycles $(123)(45)$, $G(\sigma,w)$ is the disjoint union of the previous graph and of 
$$\xymatrix{ 
(4,1)\ar[r]|1\ar[d]|2&(4,2)\ar[r]|2&(4,3)&(4,4)\ar[l]|1\\
(5,4)\ar[r]|1&(5,3)&(5,2)\ar[l]|2&(5,1)\ar[l]|1\ar[u]|2
 } 
$$


\subsubsection{Admissible graphs with restricted monochromatic  cycle and   path lengths} 

Let us fix $A_1,\ldots, A_k$ nonempty  sets of positive integers, none of them being $\{1\}$. Let  $d_1, \ldots, d_k$ denote respectively  $\sup A_1, \ldots, \sup A_k$ (which can be infinite).
  
  If an oriented edge-colored graph $G$  is admissible,   it is easy to see that for all $i\in [k]$,    the  graph obtained from $G$ by erasing all edges whose color is not $i$  is a disjoint union of    directed paths and directed cycles. If for each $i\in [k]$,  all these directed cycles  have length in $A_i$ (resp. length equal to $d_i$) and all these directed paths  have length $<d_i$, $G$ will be said to be {\it $(A_1,\ldots, A_k)$-admissible} (resp.  {\it $(d_1, \ldots, d_k)$-strongly admissible}). 
  A partition $\Delta$ of the vertex set of an oriented edge-colored graph $G$ will be said to be an {\it $(A_1,\ldots, A_k)$-admissible partition} (resp.  a {\it $(d_1,\ldots, d_k)$-strongly admissible partition}) if $G/\Delta$ is $(A_1,\ldots, A_k)$-admissible (resp.   $(d_1,\ldots, d_k)$-strongly admissible). 
  
  
  We define the {\it Neagu characteristic} of an admissible oriented edge-colored graph $G$   to be 
  $$\chi(G)=|\{\textrm{vertices of $G$}\}|-\ds\sum_{r=1}^k|\{\textrm{edges  of $G$ with color $r$}\}|+ \sum_{r=1}^k\sum_{\substack{\textrm{$L$ directed cycle}\\ \textrm{of $G$ with color $r$}}}\f{\textrm{length of $L$}}{d_r},$$ with the convention $l/\infty=0$. 

Let $G=(V;E_1,\ldots, E_k)$ be an admissible oriented edge-colored graph. A {\it  direct extension} 
of $G$ is a   graph $G'$ of one of the following types: \begin{itemize}
\item[-] $G'=(V\cup\{r\}; E_1, \ldots, E_{i-1},E_i\cup\{s\to r\},E_{i+1},\ldots, E_{k})$, with $r\notin V$,  $i\in [k]$, and $s\in V$ \st no edge of $G$ with color $i$ has $s$ for beginning vertex,
\item[-] $G'=(V\cup\{r\}; E_1, \ldots, E_{i-1},E_i\cup\{r\to s\},E_{i+1},\ldots, E_{k})$, with $r\notin V$,  $i\in [k]$, and $s\in V$ \st no edge of $G$ with color $i$ has $s$ for ending vertex,
\item[-] $G'=(V; E_1, \ldots, E_{i-1},E_i\cup\{r\to  s \},E_{i+1},\ldots, E_{k})$, with  $i\in [k]$ \st $d_i<\infty$, $r,s\in V$ \st no edge of $G$ with color $i$ has $r$ for beginning vertex, no edge of $G$ with color $i$ has $s$ for ending vertex, and there exists  $t_1,\ldots, t_{d_i}\in V$ pairwise distinct vertices \st $t_1=s, t_{d_i}=r$ and  $t_1\to t_2, t_2\to t_3,\ldots , t_{d_i-1}\to t_{d_i}\in E_i$.\end{itemize}

In other words, a direct extension of $G$ is an admissible oriented edge-colored graph which can be obtained from $G$ either by adding a vertex and   connecting this vertex to a vertex of $G$ by an edge of any sense and color, or by adding an edge which closes a monochromatic directed cycle with color $i\in [k]$ of length $d_i$. In the first case, the direct extension $G\subset G'$ is said to be {\it vertex-adding}, whereas in the second one, it is said to be {\it cycle-closing}.

For example, for $d_3=4$, the graphs $G,G', G'',G'''$ below (where the colors  of the edges appear on them) are such that $G\subset G'\subset G''\subset G'''$ are direct extensions, the first ones being vertex-adding, whereas the last one is   cycle-closing.
$$\begin{tabular}{ll}G=
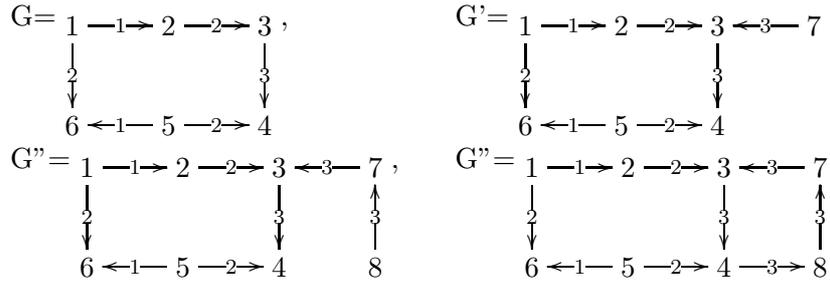
\xymatrix{1\ar[r]|1\ar[d]|2&2\ar[r]|2&3\ar[d]|3\\ 6&5\ar[l]|1\ar[r]|2&4},&\quad G'=\xymatrix{1\ar[r]|1\ar[d]|2&2\ar[r]|2&3\ar[d]|3&7\ar[l]|3\\ 6&5\ar[l]|1\ar[r]|2&4}\\ G''=\xymatrix{1\ar[r]|1\ar[d]|2&2\ar[r]|2&3\ar[d]|3&7\ar[l]|3\\ 6&5\ar[l]|1\ar[r]|2&4&8\ar[u]|3},&\quad G''=\xymatrix{1\ar[r]|1\ar[d]|2&2\ar[r]|2&3\ar[d]|3&7\ar[l]|3\\ 6&5\ar[l]|1\ar[r]|2&4\ar[r]|3&8\ar[u]|3}\end{tabular}$$


An {\it extension} of $G$ is an oriented edge-colored graph $G'$ \st there is $n\geq 0$, $G_0=G,\ldots, G_n=G'$ \st for all $i=1,\ldots,n$, $G_i$ is a direct extension of $G_{i-1}$. Such an integer $n$ is unique (it is the number of edges of $G'$ minus the number of edges of $G$) and will be called the {\it degree} of the extension.

In other words, an extension of $G$ is an admissible oriented edge-colored graph which can be obtained from $G$ by successively adding   vertices and/or edges without adding either any connected component or any monochromatic directed cycle whose length is not $d_i$ when its color is $i$.  

\begin{rmq}\label{tekilatex}It can easily be proved, using admissibility,  that a vertex or an edge of an extension of $G$ which is not a vertex of $G$ cannot belong to a monochromatic directed cycle with color $i\in [k]$ whose length is not equal to $d_i$. 
\end{rmq}

\begin{lem}\label{15.4.07.1}The Neagu characteristic is preserved   by extension.
\end{lem}

\begin{pr}It suffices to prove that the Neagu characteristic is preserved by direct extension. In the case of a vertex-adding direct extension, it is obvious. In  the case where the direct extension is monochromatic directed cycle-closing, it suffices to notice that the admissibility implies that exactly one monochromatic directed cycle is closed.
\end{pr}

The following lemma is the key result of this section. Recall that  for $u,v$ in a set $V$, $\{u=v\}$ denotes the partition of $V$ whose classes are all singletons, except one: $\{u,v\}$.
\begin{lem}\label{ultim.14.06.07}
Consider an oriented edge-colored graph $G$ which is an extension of the graph $G(w)$, with $w\in \mathbb{M}_k$ a  cyclically reduced word \st 
the order, in $F_k/[g_1^{d_1},\ldots, g_k^{d_k}]$, of the element represented by $w$ is either infinite or equal to one. Suppose that for a certain $i\in [k]$ \st $d_i<\infty$, $G$ contains an $i$-colored directed path $t_0\to t_1\to\cdots \to t_{d_i}$ with length $d_i$.  Then there is a partial $(d_1,\ldots, d_k)$-reduction $\tilde{w}$ of $w$ \st $\Adm(G/\{t_0=t_{d_i}\})$ is isomorphic to an extension of $G(\tilde{w})$. 
\end{lem}

 \begin{rmq}\label{9.06.07.1}Note that by definition of the operator $\Adm$ and by Remark  \ref{9.05.07.1},  with the notations of the previous lemma, there is a partition $\Delta$ of the vertex set of $G$ \st $\Adm(G/\{t_0=t_{d_i}\})$ can be identified with $G/\Delta$. Clearly, by definition of $\{t_0=t_{d_i}\}$ and of $\Adm$, if two vertices of $r,s$ of $G$ are  in the same class of  $\Delta$, then there are $n\geq 0$, $x_0=r, x_1, \ldots, x_{n}=s$ some vertices of $G$,  $i_1,\ldots, i_n\in [k]$, and $\eps_1,\ldots, \eps_n\in \{\pm 1\}$ such that: \\ - the $(d_1,\ldots, d_k)$-reduction of the word $g_{i_1}^{\eps_1}\cdots g_{i_n}^{\eps_n}$ is the empty word,\\ - for all $l\in [n]$, the $i_l$-colored edge $x_{l-1}\to x_l$ (resp. $x_{l-1}\leftarrow x_l$) belongs to $G$ if $\eps_l=1$ (resp. if $\eps_l=-1$). 
\end{rmq}

\begin{prlem}
 The lemma is a consequence of   the  following proposition, that we shall prove by induction on $(|w|, n)$ (for the lexical order $\preceq$): for all $w\in \mathbb{M}_k$ a cyclically reduced word \st 
the order, in $F_k/[g_1^{d_1},\ldots, g_k^{d_k}]$, of the element represented by $w$ is either infinite or equal to one,  for all extensions $G$ of degree $n$ of $G(w)$, if for a certain $i\in [k]$ \st $d_i<\infty$, $G$ contains an $i$-colored directed path $t_0\to t_1\to\cdots \to t_{d_i}$ with length $d_i$, then  $\Adm(G/\{t_0=t_{d_i}\})$ is either  isomorphic to an extension of $G({w})$ with degree $\leq n$ or isomorphic   to an extension of $G(\tilde{w})$ for a certain strict partial $(d_1,\ldots, d_k)$-reduction $\tilde{w}$ of $w$.

 For $|w|=0 $ and $n=0$, then $G=G(w)$ is a single point, thus there is nothing to prove.

 Now, let us consider a  cyclically reduced word  $w$ \st 
the order, in $F_k/[g_1^{d_1},\ldots, g_k^{d_k}]$, of the element represented by $w$ is either infinite or equal to one, and an extension $G$  of $G(w)$ with degree $n$ \st $(0,0)\prec (|w|,n)$, and let us  suppose the result to be proved for any word $w'$ and any extension $G(w')\subset G'$ with degree $n'$ \st  $(|w'|,n')\prec (|w|,n)$. Suppose that $G$ contains an $i$-colored directed path $t_0\to t_1\to\cdots \to t_{d_i}$ with length $d_i$. 


$\quad\bullet$ If   $n=0$, then $G=G(w)$, and up to a cyclic permutation of the letters of $w$ (which, up to an isomorphism, does not change $G(w)$), one can suppose that one of the following cases occurs: \begin{itemize}\item[Case 1:] $w= ug_i^{\alpha d_i}u^{-1}\tilde{w}$, with  $\alpha$ a nonzero integer, $u$ nonempty word and $\tilde{w}\in \mathbb{M}_k$ a cyclically reduced word, and there is an integer $j$, $1+|u|\leq j\leq 1+|u|+(|\alpha|-1)d_i$, \st  $$\begin{cases}t_0=j,t_1=j+1,\ldots, t_{d_i}=j+d_i&\textrm{if $\alpha >0$,}\\ t_0=j+d_i,t_1=j+d_i-1,\ldots, t_{d_i}=j&\textrm{if $\alpha <0$}.\end{cases}$$
\item[Case 2:]  $w=g_i^{\alpha d_i}\tilde{w}$, with  $\alpha$ a nonzero integer and  $\tilde{w}$ a word which can be written   $\tilde{w}=g_i^{\beta}v$, with $\beta$ an integer, $|\beta|<d_i$ and $v$ a nonempty (except possibly if $\beta=0$) word whose first and last letters do not belong to $\{g_i,g_i^{-1}\}$, and  there is an integer $j$, $1\leq j\leq 1+(|\alpha|-1)d_i+|\beta|$, \st  $$\begin{cases}t_0=j,t_1=j+1,\ldots, t_{d_i}=j+d_i&\textrm{if $\alpha >0$,}\\ t_0=j+d_i,t_1=j+d_i-1,\ldots, t_{d_i}=j&\textrm{if $\alpha <0$}.\end{cases}$$

\end{itemize}
Note that in both cases, $\tilde{w}$ is a strict partial $(d_1,\ldots, d_k)$-reduction of $w$. We are going to prove that $\Adm(G(w)/\{t_0=t_{d_i}\})$ is isomorphic to an extension of $G(\tilde{w})$.

In  Case 1, let us write $u=g_{i_1}^{\eps_1}\cdots g_{i_{|u|}}^{\eps_{|u|}}$, with $i_1,\ldots,i_{|u|}\in [k]$, $\eps_1,\ldots, \eps_{|u|}\in \{\pm 1\}$. Let us define the sequence of direct extensions 
$G(\tilde{w})=G_0\subset G_1\subset\cdots\subset G_{|u|+d_i}$ in the following way (we consider a copy $1',2',\ldots$ of the set of positive integers): 
\begin{itemize}\item[1)] First add the vertex $1'$ and the $i_1$-colored edge $1\to 1'$ or $1\leftarrow 1'$ according to whether $\eps_1=1$ or $-1$. This is   a direct extension, because since $w$ is cyclically reduced, we know that if the last (resp. first) letter of $\tilde{w}$ is $g_{i_1}^{\eps}$ with $\eps\in \{\pm 1\}$, then $\eps=\eps_1$ (resp. $-\eps_1$).
\item[2)] Then add the vertex $2'$ and the $i_2$-colored edge $1'\to 2'$ or $1'\leftarrow 2'$ according to whether $\eps_2=1$ or $-1$.  This is also a direct extension, because since $w$ is cyclically reduced, we know that if $i_1=i_2$, then $\eps_1=\eps_2$.
\item[$\vdots$]
\item[$|u|$)] Then add the vertex $|u|'$ and the $i_{|u|}$-colored edge $(|u|-1)'\to |u|'$ or $(|u|-1)'\leftarrow  |u|'$ according to whether $\eps_{|u|}=1$ or $-1$.   
\item[$|u|+1$)] Then add the vertex $(|u|+1)'$ and the $i$-colored edge $|u|'\to (|u|+1)'$.  This is   a direct extension, because we know that  $i_{|u|}\neq i$: indeed,   $w$ is cyclically reduced, and on the one hand, the letter $g_{i_{|u|}}^{\eps_{|u|}}$ appears in before $g_i^{\alpha d_i}$ in $w$, whereas on the other hand, the letter $g_{i_{|u|}}^{-\eps_{|u|}}$ appears in after $g_i^{\alpha d_i}$ in $w$.
\item[$|u|+2$)] Then add the vertex $(|u|+2)'$ and the $i$-colored edge $(|u|+1)'\to (|u|+2)'$.   
\item[$\vdots$]
\item[$|u|+d_i-1$)] Then add the vertex $(|u|+d_i-1)'$ and the $i$-colored edge $(|u|+d_i-2)'\to (|u|+d_i-1)'$.  
\item[$|u|+d_i$)] Then add the $i$-colored edge $(|u|+d_i-1)'\to |u|'$.   This is  a  direct extension because as explained at the  $(|u|+1)$th step, $i_{|u|}\neq i$.\end{itemize}
Claim: the last graph of this sequence,    $G_{|u|+d_i}$, is isomorphic to $\Adm(G(w)/\{t_0=t_{d_i}\})$. To prove this, let us first notice that by Remark \ref{9.05.07.1},  $\Adm(G/\{t_0=t_{d_i}\})$ is isomorphic to $G(w)/P$, where  $P$ is the partition of the vertex set of $G$ whose classes are singletons, except that:
\\ - In the ``$g_i^{\alpha d_i}$ part" of $G(w)$, vertices are linked every $d_i$ edges: for all $x,y\in\{|u|+1,\ldots,|u|+|\alpha |d_i+1 \}$, $x=y\mod P$ \ssi $x=y\mod d_i$.
 \\ -  In the ``$ug_i^{\alpha d_i}u^{-1}$ part" of $G(w)$, symmetric vertices of ``the $u$ and $u^{-1}$ parts"  are linked by $P$: we have, modulo $P$, $1= 2|u|+|\alpha|d_i+1$, $2=2|u|+|\alpha|d_i$,  \ldots, $|u|= |u|+|\alpha|d_i+2$. 
 
 As an illustration, on Figure \ref{0.12.09.2}, we draw the graphs $G(w)$, $G(w)/P$ and  $G_{|u|+d_i}$ for $w=g_4g_2g_3^4g_2^{-1}g_4^{-1}g_3g_1g_5g_6$,  $\tilde{w}=g_3g_1g_5g_6$ and $d_3=4$,  linking by edges of the type $\xymatrix@1{\cdot\ar@{.}[r]&\cdot}$  the vertices of $G(w)$ which are in the same class of $P$  (but these edges do not belong to  $G(w)$). \begin{figure}\begin{center}
 $$\begin{array}{ccccc}\xymatrix{&4\ar[r]|3&5\ar[r]|3&6\ar[dr]|3\\
3\ar[ur]|3\ar@{.}[rrrr]|{\textrm{linked by $P$}}&&&&7\\
2\ar[u]|2\ar@{.}[rrrr]|{\textrm{linked by $P$}}&&&&8\ar[u]|2\\
1\ar[u]|4\ar@{.}[rrrr]|{\textrm{linked by $P$}}&&&&9\ar[u]|4\ar[dl]|3\\
& 12\ar[ul]|6&11\ar[l]|5&10\ar[l]|1}&&
\xymatrix{4\ar[r]|3&5\ar[r]|3&6\ar[dl]|3\\
& \{3,7\}\ar[ul]|3\\
& \{2,8\}\ar[u]|2\\
& \{1,9\}\ar[u]|4\ar[dr]|3\\
 12\ar[ur]|6&11\ar[l]|5&10\ar[l]|1}
&&
\xymatrix{3'\ar[r]|3&4'\ar[r]|3&5'\ar[dl]|3\\
& 2'\ar[ul]|3\\
& 1'\ar[u]|2\\
& 1\ar[u]|4\ar[dr]|3\\
 4\ar[ur]|6&3\ar[l]|5&2\ar[l]|1}\\
 \textrm{The graph $G(w)$}&& \textrm{The graph $G(w)/P$}&& \textrm{The graph $G_{|u|+d_i}$}\end{array}
$$\caption{}\label{0.12.09.2}
\end{center}\end{figure} 

 Now, to be precise,  let us give a bijective map $\vfi$ from the vertex set of $G_{|u|+d_i}$ to the one of $G(w)/P$ which is an isomorphism of edge-colored oriented graphs. First, recall that the vertex set of $G_{|u|+d_i}$ is $[|\tilde{w}|]\cup\{1',\ldots, (|u|+d_i-1)'\}$ and that the vertex set of  $G(w)/P$ is  $$\{\{1, 2|u|+|\alpha|d_i+1\}, \{2, 2|u|+|\alpha|d_i\}, \{3, 2|u|+|\alpha|d_i-1\}, \ldots, \{|u|, |u|+|\alpha|d_i+2\}\}\cup$$ $$\{\{|u|+1+ ld_i\ste 0\leq l\leq |\alpha|\}\}\cup
\{\{|u|+j+ ld_i\ste 0\leq l< |\alpha|\}\ste j=2,\ldots, d_i\} \cup$$ $$\{\{2|u|+|\alpha|d_i+2\},\ldots,\{2|u|+|\alpha|d_i+|\tilde{w}| \}\}
.$$  The bijection $\vfi$ is defined by \begin{eqnarray*}\forall x\in [|\tilde{w}|],& \varphi(x)&=\begin{cases}\{1, 2|u|+|\alpha|d_i+1\}&\textrm{ if $x=1$,}\\ \{2|u|+|\alpha|d_i+x\}&\textrm{if $x>1$,}\end{cases}\\ \forall x\in [|u|+d_i-1], &\varphi(x')&=\begin{cases}\{x+1, 2|u|+|\alpha|d_i+1-x\}&\textrm{ if $x< |u|$,}\\ 
\{|u|+1+ ld_i\ste 0\leq l\leq |\alpha|\}&\textrm{if $x=|u|$,}\\
\{x+1+ld_i\ste 0\leq l<|\alpha|\}&\textrm{if $x>|u|$ and $\alpha>0$,}\\
\{2|u|+d_i+1-x+ld_i\ste 0\leq l<|\alpha|\}&\textrm{if $x>|u|$ and $\alpha<0$.}\end{cases}\end{eqnarray*}

In Case 2, 
let us define the sequence of direct extensions 
$G(v)=G_0\subset G_1\subset\cdots\subset G_{d_i-|\beta|}$ in the following way  (again,   we consider a copy $1',2',\ldots$ of the set of positive integers): 
\begin{itemize}\item[1)] First add the vertex $(|\beta|+2)'$ and the $i$-colored edge $|\beta|+1\to (|\beta|+2)'$ or $|\beta|+1\leftarrow (|\beta|+2)'$ according to whether $\beta\geq 0$ or $<0$. This is   a direct extension, because $v$ is a nonempty  (except possibly if $\beta=0$) word whose first and last letters do not belong to $\{g_i,g_i^{-1}\}$.
\item[2)] Then add the vertex $(|\beta|+3)'$ and the $i$-colored edge $(|\beta|+2)'\to (|\beta|+3)'$ or $(|\beta|+2)'\leftarrow (|\beta|+3)'$ according to whether  $\beta\geq 0$ or $<0$.  
\item[$\vdots$]
\item[$d_i-|\beta|-1$)] Then add the vertex $d_i'$ and the $i$-colored edge $(d_i-1)'\to d_i'$ or $(d_i-1)'\leftarrow  d_i'$  according to whether  $\beta\geq 0$ or $<0$.   
\item[$d_i-|\beta|$)] Then add the $i$-colored edge $d_i'\to 1$ or $d_i'\leftarrow 1$ according to whether $\beta\geq 0$ or $<0$.   \end{itemize} 

As in Case 1, one proves that  the last graph of this sequence of direct extensions, namely  $G_{d_i-|\beta|}$, is isomorphic to $\Adm(G(w)/\{t_0=t_{d_i}\})$ (which, by Remark \ref{9.05.07.1},    is isomorphic to $G(w)/P$, where  $P$ is the partition of the vertex set of $G$ whose classes are singletons, except that in the ``$g_i^{\alpha d_i+\beta}$ part" of $G$, vertices are linked every $d_i$ edges). 
As an illustration, on Figure \ref{5.12.09}, we draw the graphs $G(w)$, $G(w)/P$ and $G_{d_i-|\beta|}$, with $w=g_1^5g_2g_3g_4g_1^{-1}g_4g_{2}^{-2}$, $\tilde{w}=g_1g_2g_3g_4g_1^{-1}g_4g_{2}^{-2}$ and $d_1=4$,  linking by edges of the type $\xymatrix@1{\cdot\ar@{.}[r]&\cdot}$  the vertices of $G(w)$ which are in the same class of $P$  (but these edges do not belong to  $G(w)$). 
 \begin{figure}\begin{center}
$$\begin{array}{ccc}\xymatrix{&3\ar[r]|1&4\ar[dr]|1\\
2\ar[ur]|1\ar@{.}[drrr]&&&5\ar[d]_1\\
1\ar[u]_1\ar@{.}[urrr]|{\textrm{linked by $P$}}\ar[d]|2&&&6\ar[d]|2\\
12\ar[d]|2&&&7\ar[d]|3\\
11&&&8\ar[dl]|4\\
&10\ar[ul]|4\ar[r]|1&9
}&
\xymatrix{\\
&4\ar[d]_1&3\ar[l]|1\\
&\{1,5\}\ar[dl]|2\ar[r]|1&\{2,6\}\ar[u]_1\ar[dr]|2\\
12\ar[d]|2&&&7\ar[d]|3\\
11&&&8\ar[dl]|4\\
&10\ar[ul]|4\ar[r]|1&9}
&
\xymatrix{\\ &4'\ar[d]_1&3'\ar[l]|1\\
&1\ar[dl]|2\ar[r]|1&2\ar[u]_1\ar[dr]|2\\
8\ar[d]|2&&&3\ar[d]|3\\
7&&&4\ar[dl]|4\\
&6\ar[ul]|4\ar[r]|1&5}\\
 \textrm{The graph $G(w)$}& \textrm{The graph $G(w)/P$}&\textrm{The graph $G_{d_i-|\beta|}$}\end{array}$$\caption{}\label{5.12.09}
\end{center}\end{figure}


$\quad\bullet$ If the degree $n$ of the extension $G$ of $G(w)$ is $\geq 1$, let us introduce a sequence  $G_0=G(w)\subset G_1\subset \cdots\subset G_n= G$ of direct extensions.  Again, two cases will be to consider.

$\quad\quad\bullet\bullet$ 
  If $t_0,\ldots,t_{d_i}$ are not all vertices of $G_{n-1}$:

 Since $G_n$ is a direct extension of $G_{n-1}$, at most one vertex of $G_n$ is not a vertex of $G_{n-1}$. Hence exactly one vertex of $G_n$ is not a vertex of $G_{n-1}$,  the direct extension $G_{n-1}\subset G_n$ is vertex-adding and  the vertex added is at the extremity of exactly one edge.  Thus this vertex is either $t_0$ or $t_{d_i}$. We suppose that it is $t_{d_i}$: the case where it is $t_0$ can be treated analogously. 
 \begin{itemize}
 \item[-] If no $i$-colored edge of $G_{n-1}$ has $t_0$ for ending vertex, then $\Adm(G_n/\{t_0=t_{d_i}\})$ can simply be identified with the graph obtained from $G_{n-1}$ by adding the $i$-colored edge $t_{d_i-1}\to t_0$: this graph is a cycle-closing direct extension of $G_{n-1}$, hence an extension of degree $n$ of $G(w)$, so the result holds. 
 
\item[-] If $t_0$ is the ending vertex of a certain $i$-colored edge $t_{-1}\to t_0$ of $G_{n-1}$, then notice first that by Lemma \ref{Adm(ext)=Adm(quotient)}, $\Adm(G_n/\{t_0=t_{d_i}\})$ is isomorphic to $\Adm(G_{n-1}/\{t_{-1}=t_{d_i-1}\})$. Hence by the induction hypothesis, it suffices to prove that $t_{-1},\ldots, t_{d_i-1}$ are pairwise distinct: 
 if it wasn't the case, since $t_0, \ldots, t_{d_i-1}$ are pairwise distinct, then we would have $t_{-1}=t_l$ for some $l\in \{0,\ldots, d_i-1\}$, but it would imply that $t_l$ is the beginning vertex of two distinct $i$-colored edges of $G_n$ (the edges $t_{-1}\to t_0$ and $t_l\to t_{l+1}$), which is impossible since $G_n$, as an extension of $G(w)$, is admissible. 


  \end{itemize}
  
$\quad\quad\bullet\bullet$  If $t_0,\ldots,t_{d_i}$ are all vertices of $G_{n-1}$:
 
 There will be no ambiguity for us to also denote by $\{t_0= t_{d_i}\}$ the restriction of the partition $\{t_0=t_{d_i}\}$ to the vertex set of $G_{n-1}$. 
 By definition of the operator $\Adm$ and by Remark \ref{9.05.07.1}, there is a partition $\Gamma$ of the vertex set of $G_{n-1}$ \st $\Adm(G_{n-1}/\{t_0=t_{d_i}\})$ is isomorphic to $G_{n-1}/\Gamma$. Let, for $x$ vertex of $G_{n-1}$, $\overline{x}$ denote the class of $x$ in $\Gamma$.  We know, by the induction hypothesis, that $G_{n-1}/\Gamma$  is either  isomorphic to an extension of $G({w})$ with degree $\leq n-1$ or isomorphic   to an extension of $G(\tilde{w})$ for a certain strict partial $(d_1,\ldots, d_k)$-reduction $\tilde{w}$ of $w$.   Hence it suffices to prove one of the following properties: 
 \\
 \\
 (P1) $\Adm(G_n/\{t_0=t_{d_i}\})$ is isomorphic to  $G_{n-1}/\Gamma$, \\
 (P2)   $\Adm(G_n/\{t_0=t_{d_i}\})$ is isomorphic to a direct extension of $G_{n-1}/\Gamma$, \\
 (P3) there is $j\in [k]$ \st $d_j<\infty$ and $G_{n-1}/\Gamma$ contains a $j$-colored directed path with length $d_j$ and extremity vertices $x,y$ and \st  $\Adm(G_n/\{t_0=t_{d_i}\})$ is isomorphic to $\Adm((G_{n-1}/\Gamma)/\{x=y\})$.
 \\
 \\
 (The conclusion, in the (P3) case, requires to use the induction hypothesis again.)

 Note that  for all $\tilde{w}$  partial $(d_1,\ldots, d_k)$-cyclic reductions of $w$, since the order of the element of $F_k/[g_1^{d_1},\ldots, g_k^{d_k}]$ represented by $w$ (hence by $\tilde{w}$)  cannot be finite and nonzero, we have the following property:\\ 
(Q) for all $j\in [k]$ \st $d_j<\infty$, $G(\tilde{w})$ does not contain any $j$-colored directed cycle. 
  
Now, again, we have two sub-cases to consider:

 \begin{itemize}
 \item[-] If the extension $G_{n-1}\subset G_n$ is vertex-adding: there is $j\in [k]$, a vertex $r$ of $G_n$, a vertex $s$ of $G_{n-1}$ \st $G_n$ can be obtained from $G_{n-1}$ by the addition of the vertex $r$ and of  the $j$-colored edge $s\to r$  (or $r\to s$). Now, if no edge with color $j$ has $\overline{s}$ for beginning vertex (or respectively ending vertex) in $G_{n-1}/\Gamma$, then by Remark \ref{ext<->adm},  $\Adm(G_n/\{t_0=t_{d_i}\})$ is  isomorphic to the vertex-adding direct extension of $G_{n-1}/\Gamma$ obtained by adding the vertex $r$ and the edge $\overline{s}\to r$ (or respectively $r\to\overline{s} $). On the other hand,  if an edge with color $j$ has $\overline{s}$ for beginning vertex (or respectively ending vertex) in $G_{n-1}/\Gamma$, then $\Adm(G_n/\{t_0=t_{d_i}\})$ is   isomorphic to $G_{n-1}/\Gamma$. In both cases,  (P1) or (P2) holds.
\item[-] If the extension $G_{n-1}\subset G_n$ is   cycle-closing: there is $j\in [k]$ \st $G_n$ can be obtained from $G_{n-1}$ by adding an edge with color $j$ between two vertices of $G_{n-1}$ which are the extremities of a $j$-colored directed path with length $d_j-1$. Let us denote by $v_1,\ldots,v_{d_j}$ the successive vertices of this directed path: the edge added in the direct extension $G_{n-1}\subset G_n$ is $v_{d_j}\to v_1$, with color $j$. Notice that  $\overline{v_1},\ldots,\overline{v_{d_j}}$ are pairwise distinct vertices of  $G_{n-1}/\Gamma$. Indeed, if it wasn't the case, since $G_{n-1}$ contains the $j$-colored directed path $v_1\to\cdots\to v_{d_j}$, $G_{n-1}/\Gamma$ would contain a $j$-colored monochromatic directed cycle with length $<d_j$, which is impossible by property (Q) and Remark \ref{tekilatex}.  Again, two cases will be to consider.  \begin{itemize}
 \item[a)] If no $j$-colored edge of $G_{n-1}/\Gamma$ has $\overline{v_{d_j}}$ for beginning vertex or $\overline{v_1}$ for ending vertex: then, by Remark \ref{ext<->adm},   $\Adm(G_n/\{t_0=t_{d_i}\})$ is isomorphic to the graph obtained from $G_{n-1}/\Gamma$ by adding the vertex $\overline{v_{d_j}}\to \overline{v_1}$, which is   a   cycle-closing extension 
of $G_{n-1}/\Gamma$: (P2) holds.
 \item[b)] If a $j$-colored edge of $G_{n-1}/\Gamma$ has $\overline{v_{d_j}}$ for beginning vertex or $\overline{v_1}$ for ending vertex. We can suppose that  a $j$-colored edge of $G_{n-1}/\Gamma$ has  $\overline{v_1}$ for ending vertex: the other case can be treated analogously. Hence there is a vertex $V_{0}$ of $G_{n-1}/\Gamma$ \st $G_{n-1}/\Gamma$ contains the $j$-colored edge $V_0\to\overline{v_1}$. Now, notice that by Lemma \ref{Adm(ext)=Adm(quotient)}, $G_n/\{t_0=t_{d_i}\}$ is isomorphic to  $(G_{n-1}/\Gamma)/\{\overline{v_{d_j}}=V_0\}$: (P3) holds.
\end{itemize}
\end{itemize}
\end{prlem}


The following lemma is a consequence of the previous one, by induction on the sum, over $j\in [k]$, on   the number of $j$-colored paths of $G$ with length $\geq d_j$.

\begin{lem}\label{14.06.07.113}Consider a cyclically reduced word $w\in \mathbb{M}_k$   \st if $w$ admits a $(d_1,\ldots, d_k)$-cyclic reduction of the type $g_i^\alpha$, with $i\in [k]$ and $\alpha$ a nonzero integer, then $d_i=\infty$ and $|\alpha|\in A_i$.  Consider an oriented edge-colored graph $G$ which is an extension of the graph $G(w)$. Then $G$ admits an $(A_1,\ldots, A_k)$-admissible partition $\Delta$ \st $G/\Delta$ is isomorphic to an extension of $G(w_{\textrm{red}})$, for a certain $(d_1,\ldots, d_k)$-reduction  $w_{\textrm{red}}$ of $w$ and \st if two vertices of $r,s$ of $G$ are  in the same class of  $\Delta$, then there are $n\geq 0$,  $t_0=r, t_2, \ldots, t_{n}=s$ some vertices of $G$,  $i_1,\ldots, i_n\in [k]$, and  $\eps_1,\ldots, \eps_n\in \{\pm 1\}$ such that: \begin{itemize}\item[-] the $(d_1,\ldots, d_k)$-reduction of the word $g_{i_1}^{\eps_1}\cdots g_{i_n}^{\eps_n}$ is the empty word,\item[-] for all $l\in [n]$, the $i_l$-colored edge $t_{l-1}\to t_l$ (resp. $t_{l-1}\leftarrow t_l$) belongs to  $G$ if $\eps_l=1$ (resp. if $\eps_l=-1$).\end{itemize}
\end{lem}

The previous lemma allows us to prove  the following Theorem, which is the main result of this section.

\begin{Th}\label{margo's.mix.04.07}Consider a cyclically reduced word $w=g_{i_1}^{\alpha_1}\cdots g_{i_{|w|}}^{\alpha_{|w|}}\in \mathbb{M}_k$ \st 
 if $w$ admits a $(d_1,\ldots, d_k)$-cyclic reduction of the type $g_i^\alpha$, with $i\in [k]$ and $\alpha$ a nonzero integer, then $d_i=\infty$ and $|\alpha|\in A_i$. Then $G(w)$ admits an $(A_1,\ldots, A_k)$-admissible partition  $\Delta$ such that: \begin{itemize}
\item[(i)] $\chi(G(w)/\Delta)=\begin{cases}1&\textrm{if $w$ admits the empty word for  $(d_1,\ldots, d_k)$-cyclic reduction,}\\ 
0&\textrm{in the other case,}\\ 
\end{cases}$
\item[(ii)] for all $r<s\in [|w|]$, $r= s\mod \Delta$ implies that one of the words $$g_{i_r}^{\alpha_r}\cdots g_{i_{s-1}}^{\alpha_{s-1}}\,,\quad g_{i_s}^{\alpha_s}\cdots g_{i_{|w|}}^{\alpha_{|w|}}g_{i_1}^{\alpha_1}\cdots g_{i_{r-1}}^{\alpha_{r-1}}$$ admits the empty word as a $(d_1,\ldots, d_k)$-cyclic reduction.
\end{itemize}
\end{Th}

\begin{pr}
Note that by hypothesis, for any  $(d_1,\ldots, d_k)$-cyclic reduction $w_{\textrm{red}}$ of $w$,    $G(w_{\textrm{red}})$ is an $(A_1,\ldots, A_k)$-admissible graph with  Neagu characteristic equal to $1$ if  the only $(d_1,\ldots, d_k)$-cyclic reduction of $w$ is the empty word, and to $0$ in the other case. Hence by Lemma \ref{15.4.07.1}, for any partition $\Delta$ of the vertex set of $G(w)$, to have (i), it suffices   to prove that $G(w)/\Delta$ is an extension of $G(w_{\textrm{red}})$. Hence this Theorem is an  immediate consequence of the previous lemma, applied for $G=G(w)$.
  \end{pr}
 
  We shall use the following Corollary later. For all $l\geq 1$, $(1\cdots l)$ denotes the cyclic permutation of $[l]$ which maps $1$ to $2$, $2$ to $3$,\ldots, $l-1$ to $l$ and $l$ to $1$.
 
 \begin{cor}\label{recroisee.9.11.06}Consider a positive integer $l$ and  a cyclically reduced word $v
 \in \mathbb{M}_k$.
 
 a) Suppose that the order, in $F_k/[g_1^{d_1},\ldots, g_k^{d_k}]$, of the element represented by $v$ is  infinite and that in the case where $v$ admits a $(d_1,\ldots, d_k)$-cyclic reduction of the type $g_i^\alpha$, with $i\in [k]$ and $\alpha$ an integer, we have $l|\alpha|\in A_i$.   Then  $G((1\cdots l),v)$ admits an $(A_1,\ldots, A_k)$-admissible partition $\Delta$ such that: \begin{itemize}
\item[(i)] $\chi(G((1\cdots l),v)/\Delta)=0 $,
\item[(ii)] for all $m\neq m'\in[l]$, $
(m,1)\neq (m',1)\mod \Delta$.
\end{itemize} 

b) Suppose that the order, in $F_k/[g_1^{d_1},\ldots, g_k^{d_k}]$, of the element represented by $w$ is  equal to $l$.    Then  $G((1\cdots l),v)$ admits an $(A_1,\ldots, A_k)$-admissible partition $\Delta$ such that: \begin{itemize}\item[(i)] $\chi(G((1\cdots l),v)/\Delta)=1$,\item[(ii)] for all $m\neq m'\in[l]$, $(m,1)\neq (m',1)\mod \Delta$.\end{itemize}\end{cor}

\begin{pr} Note  that by Lemma \ref{decomp.de.base}, for all $l\geq 1$, the function $(m, i)\in [l]\tii [|v|]\mapsto (l-m)|v|+i\in [l|v|]$
 realizes an  isomorphism between $G((1\cdots l), v)$ and  $G(v^l)$, hence we are going to work with $G(v^l)$ instead of $G((1\cdots l), v)$ (and the condition ``for all $m\neq m'\in [l]$, $(m,1)\neq (m',1)\mod \Delta$" gets ``for all $m\neq m'\in [l]$, $(m-1)|v|+1\neq (m'-1)|v|+1\mod \Delta$"). 
  Then it suffices to apply Theorem \ref{margo's.mix.04.07}.
\end{pr}

\section{Words in random permutations}\label{just.for.me.12.10.06} We  fix, until the end of the article, $k\geq 1$ and $A_1$,\ldots, $A_k$ nonempty sets of positive integers, none of them being  $\{1\}$, satisfying \eqref{tech.as.14.12.07}. We shall only consider some positive integer $n$ \st  $\Sy_n(A_1)$, \ldots, $\Sy_n(A_k)$ are all  nonempty  (which is equivalent,  for $n$ large enough, to the fact that for all $i\in [k]$, $n$ is divisible by the greatest common divisor of $A_i$ \cite[Lem. 2.3]{neagu05}). 
For such an integer $n$,  we consider an independent $k$-tuple $s_1(n),\ldots, s_k(n)$ of random permutations chosen uniformly in respectively $\Sy_n(A_1)$, \ldots, $\Sy_n(A_k)$.

We also fix   a cyclically reduced word $w\in \mathbb{M}_k$ and define $\sigma_n=w(s_1(n),\ldots, s_k(n))$: the permutation obtained by replacing any $g_i$ (or $g_i^{-1}$)in $w$ by $s_i(n)$ (or $s_i(n)^{-1}$).  
 

\subsection{A key preliminary result}\label{applic.12.10.06}

 \begin{propo}\label{propo.les.pates.cuisent.et.beber.attend.10.06}Consider $p\geq 1$ and $\sigma\in \Sy_p$. The probability of the event \begin{equation}\label{event.12.10.06} \{\forall m=1,\ldots, p, \sigma_n(m)=\sigma(m)\}
 \end{equation}is equivalent, as $n$ goes to infinity,  to 
 \begin{equation}\label{eq.infini.proba}\ds\ff{n^p}\sum_{\substack{\Delta\in C(\sigma, w,A_1,\ldots, A_k)}}n^{\chi(G(\sigma, w)/\Delta)} ,\end{equation}
where $C(\sigma, w,A_1,\ldots, A_k)$ is the set of $(A_1, \ldots, A_r)$-admissible partitions $\Delta$ of $G(\sigma, w)$ such that for all $m\neq m'\in [p]$, $(m,1)\neq (m',1)\mod \Delta$.
  \end{propo}
 
 \begin{pr}
 Set  $w=g_{i_1}^{\alpha_1}\cdots g_{i_{|w|}}^{\alpha_{|w|}}$, 
 with $i_1,\ldots, i_{|w|}\in [k]$ and
  $ \alpha_1,\ldots, \alpha_{|w|}\in \{-1,1\}$. 
  We denote by $V$ the vertex set  $[p]\times [|w|]$ of $G(\sigma, w)$, fix $n\geq 1$, and define, for any 
  $s=(s_1,\ldots, s_k)\in (\Sy_n)^k$, the function 
  $$\gamma_s : (m,l)\in V\mapsto s_{i_l}^{\alpha_l}\cdots s_{i_{|w|}}^{\alpha_{|w|}}(m)\in [n].$$ Note 
  that since all $s_i$'s are one-to-one,  $\Part(\gamma_s) $  is an admissible partition of $G(\sigma, w)$. 
  Note also that if $\gamma$ is a fixed function  from  $V$ to $[n]$ \st for all $m=1,\ldots, p$, $\gamma(m,1)=\sigma(m)$, then for all $s=(s_1,\ldots, s_k)\in (\Sy_n)^k$, 
  $\gamma_s=\gamma$  \ssi  \begin{equation}\label{quick.05.1208}\textrm{for all $r\in [k]$, for all edge $e$ of $G(\sigma,w)$ with color $r$, $s_r(\gamma(\End(e)))=\gamma(\Beg(e))$.}\end{equation} Note at last that since the joint distributions of the $s_i(n)$'s is invariant under conjugation,  for such a function $\gamma$, the probability of the event $\{\gamma_{(s_1(n),\ldots, s_k(n))}=\gamma\}$ only depends on $\Part(\gamma)$.


 
 
 Hence  \pro of the event of \eqref{event.12.10.06} is the sum, over all admissible partitions $\Delta$ of $G(\sigma, w)$, of the number 
 of functions $\gamma: V\to [n]$ whose level set partition is $\Delta$  and which satisfy $ \gamma(m,1)=\sigma(m)$ for all $m\in [p]$, times the \pro that $\gamma_{(s_1(n),\ldots, s_k(n))}$ is a certain (fixed, but the choice is irrelevant) of these functions. 
 
 Now, note   that for $\Delta$  admissible partition of $G(\sigma, w)$, the number of such functions is $$\begin{cases}n(n-1)\cdots (n-|\Delta|+p+1)&\textrm{if $\forall m\neq m'\in [p], (m,1)\neq (m',1)\mod \Delta,$}\\
 0&\textrm{in the other case.}
 \end{cases}$$
 Suppose this number to be nonzero.  Then 
 by \eqref{quick.05.1208}, the \pro that $\gamma_{s_1(n),\ldots, s_k(n)}$ is a certain (fixed) of these functions is  $$\begin{cases}\prod_{r=1}^k p_n^{(A_r)}((G(\sigma,w)/\Delta)[r]) &\textrm{if $\Delta$ is an $(A_1,\ldots, A_k)$-admissible partition of $G(\sigma, w)$,}\\
  0&\textrm{in the other case,}
 \end{cases}$$
where  for all $r\in [k]$,  $(G(\sigma,w)/\Delta)[r]$ is the graph obtained from $G(\sigma,w)/\Delta$ by removing all edges which do not have color $r$ and all vertices which are not the extremity of an $r$-colored edge, and where
 for all set $A$ of positive integers and $F$ monochromatic oriented graph, $p_n^{(A)}(F)$ is the number defined in 3.8(a) of \cite{neagu05}. 
 
 So by definition of $C(\sigma, w,A_1,\ldots, A_k)$, the \pro of the event of \eqref{event.12.10.06}  is equal to $$\sum_{\substack{\Delta\in C(\sigma, w,A_1,\ldots, A_k)}}n(n-1)\cdots (n-|\Delta|+p+1)\prod_{r=1}^k p_n^{(A_r)}((G(\sigma,w)/\Delta)[r]), $$ and Proposition 3.8 of \cite{neagu05} allows one to conclude the proof. \end{pr}

The following result is a direct application of Corollary 1.3 of \cite{fbg.1}.

\begin{cor}\label{coro.ki.tue.INGdirect} If $w$ cyclically reduced is such that for all $p\geq 1$, for all $\sigma\in \Sy_p$, for all $\Delta\in C(\sigma, w,A_1,\ldots, A_k)$, one has  $\chi (G(\sigma, w)/\Delta)\leq 0$, with equality for exactly one $\Delta\in C(\sigma, w,A_1,\ldots, A_k)$, then for all $l \geq 1$, the law of $(N_1(\sigma_n),\ldots, N_l(\sigma_n))$ converges weakly, as $n$  goes to infinity,  to $$\Poiss(1/1)\otimes\cdots \otimes \Poiss(1/l).$$
\end{cor}

\begin{rmq}\label{Wolfgang+WuTang.11.06} The hypotheses of this corollary do not always hold: for $w=g_1^3g_2$, $p=1$, with $A_1=\{3,4\}$, $A_2=\{1,2\}$, the quotient of  the graph $G(w)$ by  the partition $\{1=4\}$  has Neagu characteristic $3-4+3/4+1/2>0$.\end{rmq}

\subsection{Rate of growth of the  $N_l(\sigma_n)$'s}
   \subsubsection{Case of a word with infinite order in $F_k/[g_1^{d_1},\ldots, g_k^{d_k}]$: existence of cycles of most lengths}\label{just.for.me.7.11.06}
Here, we shall see  that  even though the lengths of the cycles of $s_1(n)$,\ldots, $s_k(n)$ are supposed to belong to the specific sets $A_1,\ldots, A_k$ of positive integers, these cycles are going to mix sufficiently well to give birth to cycles of most lengths in $\sigma_n$, at least as much as in a uniform random permutation.

  \begin{Th}\label{15.11.06.989} Suppose that the order,  in $F_k/[g_1^{d_1},\ldots, g_k^{d_k}]$, of the element represented by $w$ is infinite. Then two cases can occur:
  \begin{itemize}\item[(i)]  None of the $(d_1,\ldots, d_k)$-cyclic reductions of $w$ is   of the type $g_i^\alpha$, with $i\in [k]$, $\alpha $ an integer: then for all $l\geq 1$, 
    as $n$ goes to infinity,  $$\ds\liminf( \E(N_l(\sigma_n)))\geq \ff{l}.$$
\item[(ii)]  There is $i\in [k]$ \st $A_i$ is infinite and the $(d_1,\ldots, d_k)$-cyclic reduction of $w$ is $g_i^{\alpha}$, with $\alpha$ nonzero  integer: then for all $l\geq 1$ \st $l|\alpha|\in A_i$, 
    as $n$ goes to infinity,  $$\ds\liminf (\E(N_l(\sigma_n)))\geq \ff{l}.$$
\end{itemize}
    \end{Th}  
 
 \begin{rmq}\label{bombe.de.son.15.11.06}
 In fact, we prove the more general result:  if $c_l$ denotes the cycle $(1\cdots l)$, then\begin{equation}\label{diner.alm-mcg}\forall\Delta\in C(c_l, w, A_1,\ldots, A_k),\quad\liminf \ds\lf(\f{\E(N_l(\sigma_n))}{n^{\chi(G(c_l,w)/\Delta)}}\ri)\geq \ff{l}.\end{equation} \end{rmq}

\begin{pr}
We have $$\E(N_l(\sigma_n))\!\!=\!\!\ff{l}\E\!\!\lf[\sum_{i=1}^n1_{\textrm{$i$ belongs to a cycle of length $l$ of $\sigma_n$}}\ri]\!\!=\!\!\ff{l}\sum_{i=1}^nP(\{\textrm{$i$ belongs to a cycle of length $l$}\}).$$  But in the last sum, since the law of $\sigma_n$ is invariant under conjugation, all terms are equal. Moreover, by this invariance principle again, each term is equal to the number of cycles of length $l$ containing $1$ times the probability that $\sigma_n$ contains the cycle $c_l:=(1\cdots l)$. 
So \begin{eqnarray*}\E(N_l(\sigma_n))&=&\f{n}{l}{n\choose l-1}(l-1)!P(\{\forall m=1,\ldots, l-1, \sigma_n(m)=c_l(m)\}).
 \end{eqnarray*} 
Thus, by  Proposition \ref{propo.les.pates.cuisent.et.beber.attend.10.06}, \eqref{diner.alm-mcg} follows, and then   Corollary \ref{recroisee.9.11.06} allows us to conclude.
\end{pr}

\subsubsection{Case of a word with finite order in $F_k/[g_1^{d_1},\ldots, g_k^{d_k}]$}\label{just.for.me.15.6.07}
The following theorem states that in the case where the element of  $F_k/[g_1^{d_1},\ldots, g_k^{d_k}]$  represented by $w$ has finite order $d\geq 1$,    $\sigma_n$ is not far away from having order $d$: the cardinality of the subset of $\{1,\ldots, n\}$ covered by the supports of the cycles with length $d$ in such  a random permutation  is asymptotic to $n$.
 \begin{Th}\label{15.6.07.1} Suppose that  the element of  $F_k/[g_1^{d_1},\ldots, g_k^{d_k}]$  represented by $w$ has finite order $d\geq 1$. Then,      as $n$ goes to infinity, $\Sy_n(A_r)\neq \emptyset$,   for all $p\in [1, +\infty)$, $N_d(\sigma_n)/n$ converges to $1/d$ in $L^p$ and for all $l\neq d$, $N_{l}(\sigma_n)/n$ converges to $0$ in $L^p$.
   \end{Th}  

\begin{pr}Notice first that   for all $n$, \begin{equation}\label{diner.alm-mcg.2}\sum_{l\geq 1}l N_{l}(\sigma_n)/n= 1, \end{equation}which implies that  for all positive integer $l$, the sequence $N_{l}(\sigma_n)/n$ is bounded by $0$ and $1/l$. Now, recall that for a sequence $(X_n)$ of random variables, if there exists $M<\infty$ \st for all $n$, $0\leq X_n\leq M$, the convergence of $X_n$ to a limit $X$ in all $L^p$ spaces is implied by this result for $p=1$ (indeed, under this hypothesis, if for a certain $p\geq 1$, $X_n$ does not converge to $X$ in $L^p$, then a subsequence of $(X_{n})$ converges to $X$ almost surely but not in $L^p$, which is impossible by the Dominated Convergence Theorem).   

Hence by \eqref{diner.alm-mcg.2}, it suffices to prove that $\liminf (\E(N_d(\sigma_n)/n))\geq 1/d$. This follows from   \eqref{diner.alm-mcg} (which has been established in the proof of Theorem \ref{15.11.06.989} without using the specific hypothesis of this theorem) and of Corollary \ref{recroisee.9.11.06}. \end{pr}

\subsection{Case when all $A_i$'s are infinite}\label{5.2.07.1}
 
 \begin{Th}\label{2046.15.11.06.1} Suppose that all 
 $A_i$'s are infinite, that $|w|>1$ and that $w$ is not a power of another word. Then as $n$  goes to infinity, for all $l\geq 1$,  the law of $(N_1(\sigma_n),\ldots, N_l(\sigma_n))$ converges weakly to $$\Poiss(1/1)\otimes\cdots \otimes \Poiss(1/l).$$
 \end{Th}
  
\begin{pr} This result is a direct consequence of Corollary \ref{coro.ki.tue.INGdirect}. Note first that since all $A_i$'s are infinite, for all $i$, $d_i=\infty$. Hence for all $p\geq 1$, for all $\sigma\in \Sy_p$, for all $\Delta$ an admissible partition of  $G(\sigma, w)$,   $\chi (G(\sigma, w)/\Delta)$ is the number of classes of 
$\Delta$ minus the number of edges of $G(\sigma, w)/\Delta$, hence is not positive by  Proposition 2.4.3  of \cite{nica94}. The   last sentence of Corollary 2.6.6 of the same article   also says that there is only one admissible partition $\Delta$ of $G(\sigma, w)$   \st  $\chi (G(\sigma, w)/\Delta)$  is null and \st  for all $m\neq m'\in [p]$, $(m,1)$ and $(m',1)$ are not in the same class: it is the singletons partition. It remains only to prove that the singletons partition is in $C(\sigma, w,A_1,\ldots, A_k)$:

- this partition is  admissible and does not link $(m,1)$ and $(m',1)$ for $m\neq m'$,

- there is no monochromatic directed cycle in $G(\sigma, w)$: indeed,  by Lemma \ref{decomp.de.base}, $G(\sigma, w)$ is a disjoint union of graphs of the type $G(w^d)$ ($d\geq 1$), where there is no monochromatic directed cycle since $w$ is primitive and $|w|>1$.
\end{pr}

The previous theorem, with the fact that 
   for any permutation $\sigma$ and $l,p\geq 1$, we have $N_l(\sigma^p)=\ds\sum N_{lp/h}(\sigma)p/h,$ where the sum runs over the positive integers $h$ \st $h|p$ and $\operatorname{gcd}(h,l)=1$, allows us to compute easily the limit distribution of $(N_1(\sigma_n),\ldots, N_l(\sigma_n))$ whenever $w$ is not a power of a $g_i$.  
  
\subsection{Case where $w=g_1\cdots g_k$} Now, we are not going to make the hypothesis  that all $A_i$'s are infinite anymore, but we are going to suppose that $w$ is a particular word: $w=g_1\cdots g_k$, with $k\geq 2$ (the case $k=1$ has already been treated in \cite{fbg.1}). 
\subsubsection{Case where $k>2$ or $A_1\cup A_2\nsubseteq\{1,2\}$}

\begin{Th}\label{downpressure.peter.teuchi} Under this hypothesis, as $n$  goes to infinity, for all $l\geq 1$,  the law of $(N_1(\sigma_n),\ldots, N_l(\sigma_n))$ converges weakly to $$\Poiss(1/1)\otimes\cdots \otimes \Poiss(1/l).$$
\end{Th}

 In order to prove the Theorem, we shall need the following lemmas.
 
 \begin{lem}Let $\mc{X}$ be a finite set, let $\mc{B}$ be a set of subsets of $\mc{X}$ which have all cardinality $2$. Let $\Delta$ be a partition of $\mc{X}$ such that for all $\{x,y\}\in \mc{B}$, $x=y\mod \Delta$. Then  $|\Delta|\leq |\mc{X}|-|\mc{B}|.$
\end{lem}

\begin{pr}Let us define the set $\Gamma$ of subsets of $\mc{X}$ 
by $\Gamma=\mc{B}\cup\{\{z\}\ste z\in \mc{X}\ds, z\notin 
\cup_{\{x,y\}\in \mc{B}} \{x,y\}\}.$
 By hypothesis, any class of $\Delta$ is a union of elements of $\Gamma$, so $|\Delta|\leq |\Gamma|=|\mc{X}|-|\mc{B}|.$ \end{pr}

\begin{lem}\label{SufjanStevens} Consider $w=g_1\cdots g_k$, with $k\geq 2$, and $\sigma\in \Sy_p$, with $p\geq 1$. Let $\Delta $ be an admissible partition of   $G(\sigma,w)$ 
such that for all $i\neq j\in [p]$, $(i,1)\neq (j,1)\mod \Delta$. Then \begin{itemize}\item[(i)] two different edges of $G(\sigma, w)$ with the same color cannot have their  beginning vertices in the same class of $\Delta$,  and the same holds for ending vertices,  \item[(ii)] the following inequalities hold: 
\begin{itemize}\item[(a)] if $k>2$, then $|\Delta|\leq pk-\ds\sum_{r=1}^k\sum_{\substack{\textrm{$L$ directed cycle of}\\ G(\sigma,w)/\Delta\textrm{ with color $r$}}}\textrm{length of $L$} ,$
\item[(b)] if $k=2$, then for all $r=1,2$,  $|\Delta|\leq pk-\ds\!\!\!\!\!\!\!\!\sum_{\substack{\textrm{$L$   directed cycle of}\\ G(\sigma,w)/\Delta \textrm{ with color $r$}}}\textrm{length of $L$}.$ 
\end{itemize}
\end{itemize}
\end{lem}

\begin{rmq}\label{laphilodulang.5.11.06}Note that (i) implies that there is a canonical identification between the set of edges of the graph $G(\sigma,w)$ and the set of edges of the graph $G/\Delta$: if, for all vertex $x$ of $G(\sigma,w)$, one denotes by $\overline{x}$ the class of $x$ in $\Delta$, then for all $r\in [k]$, the function from the set of $r$-colored edges of $G(\sigma,w)$ to the set of $r$-colored edges of $G(\sigma,w)/\Delta$ which maps any edge $x\to y$ to $\overline{x}\to\overline{y}$ is bijective.\end{rmq}

 \begin{pr}
 Note first that since for all $i\neq j\in [p]$, $(i,1)\neq (j,1)\mod \Delta$, and by definition of an admissible partition, one has (by an obvious induction on $l$):\begin{equation}\label{18.08.06.1} \forall i\neq j\in [p], \forall l\in [k], \;(i,l)\neq (j,l)\mod \Delta.\end{equation} 
 Since two edges of  $G(\sigma,w)$ with the same color have beginning vertices (and ending vertices) with the same second coordinate, it implies (i).
 
Now, note that, as observed in  Remark \ref{laphilodulang.5.11.06}, the edges of $G(\sigma,w)/\Delta$ can be identified with the ones of $G(\sigma,w)$. 
 Let, for $r\in[k]$, $\mc{L}[r]$ be the set of edges of $G(\sigma,w)$ which, via this identification,  belong to an $r$-colored    directed cycle of  $G(\sigma,w)/\Delta$ and  define $\ds\mc{L}:=\cup_{r=1}^k\mc{L}[r].$
 
 {\it \underline{Claim}: } no edge of $G(\sigma,w)/\Delta$ can belong to more than one monochromatic directed cycle of $G(\sigma,w)/\Delta$. 
 
 Indeed, if it where the case, there would be   $r\in [k]$,   $e,e',e''\in \mc{L}[r]$   such that $e'$ and $e''$ both follow $e$ in $r$-colored    directed cycles of $G(\sigma,w)/\Delta$ and $e'\neq e''$. Since $e'$ and $e''$ follow $e$ in   directed cycles of  $G(\sigma,w)/\Delta$, we have $\Beg(e')=\End(e)=\Beg(e'')\mod \Delta,$ which implies $e'=e''$ by (i). Contradiction.

   The first consequences of this claim are that \begin{equation}\label{11.05.07.11}\ds|\mc{L}|=\sum_{r=1}^k\sum_{\substack{\textrm{$L$   directed cycle of}\\ G(\sigma,w)/\Delta \textrm{ with color $r$}}}\textrm{length of $L$} ,\end{equation} and that for all $r=1,\ldots, k$, \begin{equation}\label{11.05.07.22}\ds |\mc{L}[r]|=\sum_{\substack{\textrm{$L$   directed cycle of}\\ G(\sigma,w)/\Delta\textrm{ with color $r$}}}\textrm{length of $L$}.\end{equation}

Another consequence of the claim is that one can define a permutation  $\varphi$  of $\mc{L}$ which maps any edge $e\in \mc{L}$ to the edge which follows $e$ in the monochromatic directed cycle  $e$ belongs to.     Let us define, for $e\in \mc{L}$, $$S(e):=\{\End(e), \Beg(\varphi(e))\}\subset [p]\tii [|w|].$$
Note that since $k\geq 2$, for all edge $e$ of $ G(\sigma,w)$, the color of $e$ is not the same as has the one of the edge whose beginning is the end of $e$. This  allows us to claim that for all $e\in \mc{L}$, $|S(e)|=2$.

For all $e\in \mc{L}$, since $\varphi(e)$ follows $e$ in $G(\sigma,w)/\Delta$, we have  $\End(e)= \Beg(\varphi(e))\mod \Delta$. So, in order to apply the previous lemma, we have to minor the cardinality of $$\mc{A}:=\{S(e)\ste e\in \mc{L}\}\quad\textrm{ (to prove (a))}$$ or, for $r\in [k]$,  of  $$\mc{A}[r]:=\{S(e)\ste e\in \mc{L},\textrm{ $e$ has color $r$}\}\quad\textrm{ (to prove (b))}.$$Suppose that there are  $e\neq f\in \mc{L}$ \st $S(e)=S(f)$. One has either $$(\End(e), \Beg(\varphi(e)))=(\End(f), \Beg(\varphi(f)))$$ or $$(\End(e), \Beg(\varphi(e)))=(\Beg(\varphi(f)),\End(f)).$$ But $(\End(e), \Beg(\varphi(e)))=(\End(f), \Beg(\varphi(f)))$ is impossible because two different edges of $ G(\sigma,w)$ cannot have the same end, since no letter of $w$ has the exponent $-1$. So one has $(\End(e), \Beg(\varphi(e)))=(\Beg(\varphi(f)),\End(f))$.

 {\it \underline{If $k>2$}: }let us prove that $S(e)=S(f)$ with $e\neq f$ is impossible. We have $\End(e)=\Beg(\varphi(f))$, so the color following the one of $e$ in the cyclic order $1,2,\ldots, k,1,\ldots$ is the one of $\varphi(f)$, i.e. of $f$. In the same way, the relation $\End(f)=\Beg(\varphi(e))$ implies that the color following the one of $f$ in the same cyclic order is the one of $e$. To sum up, in this cyclic order, one has the following direct sequence: $$\ldots,\textrm{color of $e$}, \textrm{color of $f$},\textrm{color of $e$},\ldots$$ which is impossible, since $k>2$.

So the cardinality of $\mc{A}$ is the one of $\mc{L}$, and by \eqref{11.05.07.11}, the result (a) of the lemma is an immediate application of the previous lemma, for $\mc{X}=[p]\tii [|w|]$ and $\mc{B}=\mc{A}$.

 {\it \underline{If $k=2$}: }  $\End(f)=\Beg(\varphi(e))$  implies that the color of $f$ is different from the color of $\varphi(e)$, i.e. of $e$. So for all $r=1,2$, the cardinality of  $\mc{A}[r]$ is the one of $\mc{L}[r]$, and  by \eqref{11.05.07.22}, the result (b) is an immediate application of the previous lemma, for $\mc{X}=[p]\tii [|w|]$ and $\mc{B}=\mc{A}[r]$.
 \end{pr}

 \noindent {\bf Proof of Theorem  \ref{downpressure.peter.teuchi}. }
 Again, we are going to apply  Corollary \ref{coro.ki.tue.INGdirect}. Let us fix $p\geq 1$ and $\sigma\in \Sy_p$.
Note that the singletons partition, denoted by $\Delta_s$, is in $C(\sigma, w,A_1,\ldots, A_k)$ (the proof is the same one as in the proof of Theorem \ref{2046.15.11.06.1}). Moreover, using  Lemma \ref{decomp.de.base}, 
one easily sees that the Neagu characteristic of $G(\sigma,w)$   is $0$.
Hence it suffices to prove that for all $\Delta\in C(\sigma, w,A_1,\ldots, A_k)$ \st $\Delta \neq \Delta_s$,  the Neagu characteristic of $G(\sigma, w)/\Delta$ is negative. Let us fix such a partition $\Delta$.
By Remark \ref{laphilodulang.5.11.06}, we have to prove that   \begin{equation}\label{21.08.06.12}|\Delta|< pk-\ds\sum_{r=1}^k\sum_{\substack{\textrm{$L$   directed cycle of}\\ G(\sigma,w)/\Delta \textrm{ with color $r$}}}  \f{\textrm{length of $L$}}{d_r }.
\end{equation}

$\bullet$ If there is no monochromatic directed cycle in $G(\sigma, w)/\Delta$, then, since $\Delta\neq \Delta_s$,  $|\Delta|<|[p]\tii [|w|]|=pk$, and hence (\ref{21.08.06.12}) holds.

$\bullet$ If  there is a monochromatic directed cycle in $G(\sigma, w)/\Delta$ and $k>2$, then  since for all $r$, $d_r\geq 2>1$, by  (ii) (a) of Lemma \ref{SufjanStevens}, (\ref{21.08.06.12}) holds.

$\bullet$ If  there is a monochromatic directed cycle in $G(\sigma, w)/\Delta$ and $k=2$. First note that adding  (ii) (b) of Lemma \ref{SufjanStevens}  for $r=1$ and $r=2$, and then dividing by $2$, one gets  \begin{equation}\label{21.08.06.13}|\Delta|\leq pk-\ff{2}\ds\sum_{r=1}^2\sum_{\substack{\textrm{$L$   directed cycle of}\\ G(\sigma,w)/\Delta \textrm{ with color $r$}}} \textrm{length of $L$} .
\end{equation}By hypothesis, one has either $d_1>2$ or $d_2>2$. By symetry, we will suppose that $d_1>2$: 
\begin{itemize}
\item[-]If there is at least one   directed cycle of color $1$  in $G(\sigma, w)/\Delta$. Then (\ref{21.08.06.12}) holds, because since $d_1>2$, the right hand term of (\ref{21.08.06.13}) is strictly less than $$pk-\ds\sum_{r=1}^2\sum_{\substack{\textrm{$L$   directed cycle of}\\ G(\sigma,w)/\Delta\textrm{ with color $r$}}}\f{\textrm{length of $L$}}{d_r } .$$
\item[-]If there is no   directed cycle of color $1$  in $G(\sigma, w)/\Delta$. Then there is at least one   directed cycle of color $2$, and, since $d_2>1$, 
the right hand term in (ii) (b) of Lemma \ref{SufjanStevens}  is strictly less than $$pk-\ds\sum_{\substack{\textrm{$L$   directed cycle of}\\ G(\sigma,w)/\Delta \textrm{ with color $2$}}}\f{\textrm{length of $L$}}{d_2} =pk-\ds\sum_{r=1}^2\sum_{\substack{\textrm{$L$   directed cycle of}\\ G(\sigma,w)/\Delta\textrm{ with color $r$}}}\f{\textrm{length of $L$}}{d_r},$$ i.e. than the right hand term of  (\ref{21.08.06.12}). So (\ref{21.08.06.12}) holds.
\end{itemize}{\hfill  $\square$}

 \subsubsection{Case where $k=2$ and $A_1\cup A_2\subset \{1,2\}$}Here, for $n$ positive integer, we consider $\sigma_n=s_1(n)s_2(n)$ with $s_1(n),s_2(n)$ independent random involutions. There are three cases:\begin{itemize}\item[(i)] $s_1(n),s_2(n)$ have both uniform distribution on the set of involutions of $[n]$:  $A_1= A_2=\{1,2\}$,
\item[(ii)] $s_1(n),s_2(n)$ have both uniform distribution on the set of involutions of $[n]$ without any fixed point:  $A_1= A_2=\{2\}$,
\item[(iii)] one of the random involutions has uniform distribution on the set of involutions of $[n]$ and the other one has  uniform distribution on the set of involutions of $[n]$ without any fixed point: since for $s,s'$ involutions, $ss'$ and $s's$ are conjugate,  one can suppose that $A_1=\{2\}$, $A_2=\{1,2\}$.
\end{itemize}

 \begin{Th}\label{mardi.28.11.06.2}
 \begin{itemize}\item[(i)] If $A_1= A_2=\{1,2\}$, then as $n$ goes to infinity, for all $q\geq 1$,  the law of $(N_{1}(\sigma_n), \ldots, N_{q}(\sigma_n))$ converges weakly to $\mu_1\otimes\cdots \otimes \mu_q,$ where for all $l\geq 1$, $\mu_l$ is the law of $\Po_1+2\Po_{\ff{2l}}$, with $\Po_1,\Po_{\ff{2l}}$ independent random variables with  distributions  $\Poiss(1), \Poiss(\ff{2l})$.
\item[(ii)] If  $A_1= A_2=\{2\}$, then as $n$   goes to infinity, for all $q\geq 1$,  the law of $$\ff{2}(N_{1}(\sigma_n), \ldots, N_{q}(\sigma_n))$$ converges weakly to $\Poiss(1/2)\otimes\Poiss(1/4)\otimes\cdots \otimes \Poiss(1/2q).$
\item[(iii)] If $A_1=\{2\}$, $A_2=\{1,2\}$, then  as $n$  goes to infinity, for all $q\geq 1$,  the law of $(N_{1}(\sigma_n), \ldots, N_{q}(\sigma_n))$ converges weakly to $\mu_1\otimes\cdots \otimes \mu_q,$ where for all $l$ odd,  $\mu_l$ is the law of $2\Po_{\ff{2l}}$  and for all $l$ even, $\mu_l$  
is the law of $\Po_\ff{2}+2\Po_{\ff{2l}}$, with $\Po_\ff{2},\Po_{\ff{2l}}$ independent random variables with  distributions  $\Poiss(\ff{2}), \Poiss(\ff{2l})$.
\end{itemize}
\end{Th}

To prove the Theorem, we shall need the following lemmas. 

\begin{lem}\label{1.12.06.apres.tango}Consider $a,b$ positive real numbers and define the measure $$\nu_{a,b}:=\ds e^{-\f{1+2ab}{2b^2}}\sum_{r=0}^{+\infty}\f{\E[(X+a)^r]}{r! b^r}\delta_r,$$with $X$ standard Gaussian random variable. 
Then $\nu_{a,b}$ is the distribution of $\Po_{\f{a}{b}}+2\Po_{\ff{2b^2}}$, with  $\Po_{\f{a}{b}},\Po_{\ff{2b^2}}$ independent random variables with respective distributions $\Poiss(a/b), \Poiss(1/(2b^2))$.
\end{lem}


\begin{lem}\label{mardi.28.11.06.1}Consider $A_1,A_2$ as in (i), (ii) or (iii) of Theorem \ref{mardi.28.11.06.2}.  Fix $p\geq 1$ and $\sigma\in\Sy_p$. Then the following holds: 

$\bullet$ for all  $\Delta\in C(\sigma, g_1g_2, A_1, A_2)$, $\chi(G(\sigma, g_1g_2)/\Delta)=0$, 

$\bullet$ the cardinality of  $C(\sigma, g_1g_2, A_1, A_2)$ is $$\begin{cases}\prod_{l=1}^p \E((\sqrt{l}X+l+1)^{N_l(\sigma)})&\textrm{if $A_1= A_2=\{1,2\}$,}\\
\prod_{l=1}^p \E((1+\sqrt{l}X)^{N_l(\sigma)})&\textrm{if $A_1= A_2=\{2\}$,}\\
\ds\prod_{\substack{1\leq l \leq p\\ \textrm{$l$ odd}}} \E[(\sqrt{l}X+1)^{N_l(\sigma)}]\prod_{\substack{1\leq l \leq p\\ \textrm{$l$ even}}} \E[(\sqrt{l}X+l/2+1)^{N_l(\sigma)}]&\textrm{if $A_1=\{2\}$, $A_2=\{1,2\}$,} \end{cases}$$where $X$ is a standard Gaussian random variable.
\end{lem}

\begin{pr} 
 Let, for $l$ positive integer,  $c_{2l}$ be the   cycle $(1\, 2\cdots 2l)$.
   
First, we are going to use Lemma \ref{decomp.de.base}: up to a renaming of its vertices and edges, $G(\sigma, g_1g_2)$ is the disjoint union, for $(l,i)$ \st $l\in [p], i\in [N_l(\sigma)]$, of the graphs $H(l,i)$, where each $H(l,i)$ is the oriented edge-colored graph with color set $\{1,2\}$, with vertex set $\{(l,i,j)\ste j\in [2l]\}$, with $1$-colored edge set $\{(l,i, j)\to(l,i, c_{2l}(j))\ste j\in [2l]\textrm{ odd}\}$ and with $2$-colored edge set $\{(l,i, j)\to(l,i,c_{2l}(j))\ste j\in [2l]\textrm{ even}\}$. An illustration is given in Figure \ref{11.12.08}.   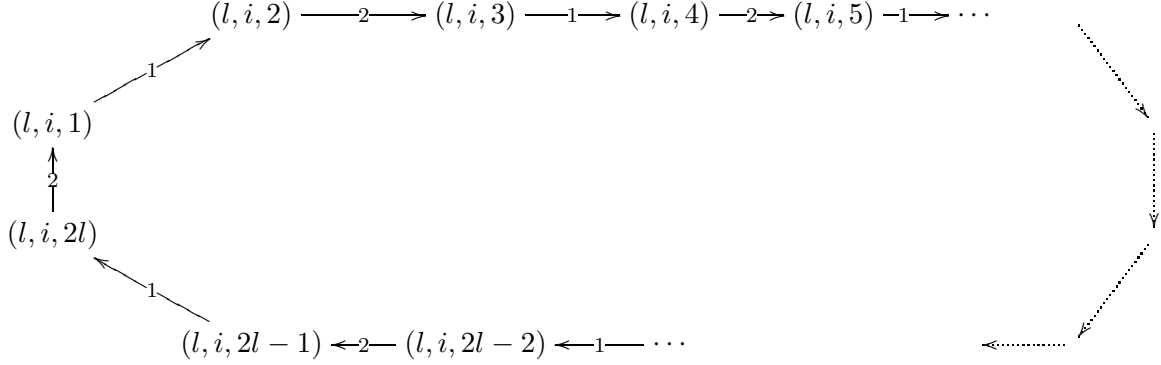
\begin{figure}\begin{center}$$\xymatrix{
&(l,i,2)\ar[r]|2&(l,i,3)\ar[r]|1&(l,i,4)\ar[r]|2&(l,i,5)\ar[r]|1&\cdots&\ar@{.>}[dr]\\
(l,i,1)\ar[ur]|1&&&&&&&\ar@{.>}[d]\\
(l,i,2l)\ar[u]|2&&&&&&&\ar@{.>}[dl]\\
&(l,i,2l-1)\ar[ul]|1&(l,i,2l-2)\ar[l]|2&\cdots\ar[l]|-{1}&&&\ar@{.>}[l]
}$$\caption{The graph $H(l,i)$.}\label{11.12.08}\end{center}\end{figure}
 With this renaming, for $\Delta$ a partition of the vertex set of $G(\sigma, g_1g_2)$, the condition $$\forall m\neq m'\in [p], (m,1)\neq (m',1)\mod \Delta$$ gets $$\forall l,l'\in [p], \forall (i,i')\in [N_l(\sigma)]\tii [N_{l'}(\sigma)], \forall (j,j')\in [2l]\tii [2l'], $$ \begin{equation}\label{27.11.06.1} [(l,i,j)\neq (l',i',j')\textrm{ and } j=j'\mod 2]\; \Rightarrow \, [(l,i,j)\neq (l',i',j')\mod \Delta].\end{equation}

Now, let us fix $\Delta\in C(\sigma, g_1g_2, A_1, A_2)$. 

a) By Remark \ref{laphilodulang.5.11.06},  $G(\sigma,w)/\Delta$ has $2p$ edges.

b) By (\ref{27.11.06.1}), the classes of $\Delta$ cannot have more than two elements. 

c) $\Delta$ is an $(A_1, A_2)$-admissible partition, so no directed path of color either $1$ or $2$ of $G(\sigma, g_1g_2)/\Delta$ can have length at least two. This implies that if two distinct  edges of $G(\sigma, g_1g_2)/\Delta$ have the same color and an extremity in common, then these edges are the ones of a directed cycle of length two, i.e. that the beginning of each of them is the end of the other one and vice versa.

d) Let us consider two distinct vertices $(l,i,j)\neq (l',i',j')$ of $G(\sigma, g_1g_2)$ which form together a class of  $\Delta$.  

\indent 
\indent ($\alpha$) By (\ref{27.11.06.1}), $j$ and $j'$ do not have the same parity, hence the edge which begins at $(l,i,j)$ has the color of the edge which ends at $(l',i',j')$, and vice versa.  So  by  an immediate induction on $|m|$ (using c)), one has:
 \begin{equation}\label{27.11.06.2} \forall m\textrm{ positive or negative integer}, \,(l,i,c_{2l}^m(j))= (l',i',c_{2l'}^{-m}(j'))\mod \Delta .\end{equation}
 
\indent \indent   ($\beta$) Suppose that $l\neq l'$. One can suppose that $l<l'$. Then by (\ref{27.11.06.2}) for $m=2l$, $(l,i,j)=(l,i,c_{2l}^{2l}(j))$ is in the same time linked, by $\Delta$, with  $(l',i',c_{2l'}^{-2l}(j'))$ and $ (l',i',j')$ which are not the same. By b), it is impossible. Hence $l=l'$. 

\indent \indent  ($\gamma_1$)Suppose that $i\neq i'$. Then  by (\ref{27.11.06.2}) for  $m=-j+1,-j+2,\ldots, -j+2l$, there is  a unique $q\in [l]$ \st the classes, in  $\Delta$, of  $$(l,i, 1), (l,i,2),\ldots, (l,i,2l) ,(l,i', 1), (l,i',2),\ldots, (l,i',2l)  $$ are\begin{equation}\label{26.06.07.1} \{(l,i, 1), (l,i',2q)\} ,\; \;\{(l,i, 2), (l,i',c_{2l}^{-1}(2q))\},\; \;\ldots, \; \;\{(l,i, 2l), (l,i',c_{2l}^{-(2l-1)}(2q))\}.\end{equation} Note that $q$ depends   on $i$ and $i'$ in a symmetric way, since the set of \eqref{26.06.07.1} can also be written $$\{(l,i', 1), (l,i,2q)\} ,\; \; \{(l,i', 2), (l,i,c_{2l}^{-1}(2q))\},\; \;\ldots, \,\{(l,i', 2l)\, ,\; \; (l,i,c_{2l}^{-(2l-1)}(2q))\}.$$
 On Figure \ref{7.12.08.1}, we give an illustration of the case $(\gamma_1)$ with $l=4$, $q=2$. In the inner (resp. outer) octagon, the vertices where not denoted by $(l,i,j)$ (resp. $(l,i',j)$) but by $j$ (resp. $j'$), for $j\in [8]$, in order to lighten the figure. In $G(\sigma, g_1g_2)$, vertices linked by $\Delta$ are linked by edges of the type $\xymatrix@1{\cdot\ar@{.}[r]|\Delta&\cdot}$ (but these edges do not belong to  $G(\sigma, g_1g_2)$). On this figure,  it appears clearly that if two vertices of the left graph linked  by an edge of the type $\xymatrix@1{\cdot\ar@{.}[r]|\Delta&\cdot}$ are in the same class of $\Delta$, then, since  no directed path of color either $1$ or $2$ of $G(\sigma, g_1g_2)/\Delta$ can have length at least two, all the other pairs of vertices  linked  by an edge of the type $\xymatrix@1{\cdot\ar@{.}[r]|\Delta&\cdot}$ must form classes of $\Delta$.
\begin{figure}\begin{center}$$\begin{array}{ccc}\xymatrix{
&&3'\ar[ddll]|-1\ar@{.}[d]|\Delta&2'\ar@{.}[d]|\Delta\ar[l]|-2\\
&&2\ar[r]|-2&3\ar[dr]|-1\\
4'\ar@{.}[r]|\Delta\ar[d]|-2&1\ar[ur]|-1&&&4\ar@{.}[r]|\Delta\ar[d]|-2&1'\ar[uull]|-1\\
5'\ar@{.}[r]|\Delta\ar[ddrr]|-1&8\ar[u]|-2&&&5\ar@{.}[r]|\Delta\ar[dl]|-1&8'\ar[u]|-2\\
&&7\ar@{.}[d]|\Delta\ar[ul]|-1&6\ar@{.}[d]|\Delta\ar[l]|-2\\
&&6'\ar[r]|-2&7'\ar[uurr]|-1
}&&\xymatrix{\\ &\{2,3'\}\ar@/_/[r]|2\ar@/_/[dl]|1&\{3,2'\}\ar@/_/[l]|2\ar@/_/[dr]|1\\
\{1,4'\}\ar@/_/[ur]|1\ar@/_/[d]|2&&&\{4,1'\}\ar@/_/[ul]|1\ar@/_/[d]|2\\
\{8,5'\}\ar@/_/[dr]|1\ar@/_/[u]|2&&&\{5,8'\}\ar@/_/[dl]|1\ar@/_/[u]|2\\
&\{7,6'\}\ar@/_/[r]|2\ar@/_/[ul]|1&\{6,7'\}\ar@/_/[l]|2\ar@/_/[ur]|1
}
\\
\substack{\textrm{restriction of the graph $G(\sigma, g_1g_2)$ to}\\
\{(l,i, 1), (l,i,2),\ldots, (l,i,8) \}\cup\{(l,i', 1), (l,i',2),\ldots, (l,i',8) \}}&&
\substack{\textrm{restriction of the graph $G(\sigma, g_1g_2)/\Delta$ to the set}\\ \textrm{of classes of }
(l,i, 1), (l,i,2),\ldots, (l,i,8) ,(l,i', 1), (l,i',2),\ldots, (l,i',8) }
\end{array}$$\caption{}\label{7.12.08.1}\end{center}\end{figure}

 \indent  \indent  ($\gamma_2$) Let us prove that if $i=i'$, then there is $j_0\in [2l]$ \st  the classes, in $\Delta$, of the elements  $(l,i, 1),  (l,i,2),\ldots, (l,i,2l)  $ are\begin{eqnarray}& \{(l,i,j_0),(l,i,c_{2l}(j_0)) \}, \;\;\{(l,i, c_{2l}^{-1}(j_0)), (l,i,c_{2l}^2(j_0))\},\;\; \{(l,i, c_{2l}^{-2}(j_0)), (l,i,c_{2l}^3(j_0))\},&\label{knives.out.30.11.06}\\ &\ldots, \;\;  \{(l,i, c_{2l}^{-(l-1)}(j_0)), (l,i,c_{2l}^l(j_0))\},&\nonumber \end{eqnarray}
 and that this $j_0$ is unique up to a replacement by $c_{2l}^l(j_0)$.
 
  Note first that if such a $j_0$ exists, then its  uniqueness up to a replacement by $c_{2l}^l(j_0)$ is obvious. Indeed, if one denotes $c_{2l}^l(j_0)$ by $j_0'$, then since $c_{2l}$ is a cycle of length $2l$, the partition of (\ref{knives.out.30.11.06}) is equal to
\begin{eqnarray*}& \{(l,i,j_0'),(l,i,c_{2l}(j'_0)) \}, \;\;\{(l,i, c_{2l}^{-1}(j'_0)), (l,i,c_{2l}^2(j'_0))\},\;\;\{(l,i, c_{2l}^{-2}(j'_0)), (l,i,c_{2l}^3(j'_0))\},&\\ &\ldots,\;\;  \{(l,i, c_{2l}^{-(l-1)}(j'_0)), (l,i,c_{2l}^l(j'_0))\},&\end{eqnarray*} and $\{(l,i,j_0),(l,i,c_{2l}(j_0)) \}$,  $ \{(l,i, c_{2l}^{-(l-1}(j_0)), (l,i,c_{2l}^l(j_0))\}$ ($=\{(l,i,j'_0)$, $(l,i,c_{2l}(j'_0)) \}$) are the only classes of the type $\{(l,i,x),(l,i,c_{2l}(x)) \}$ in the partition of (\ref{knives.out.30.11.06}).  

To prove the existence of such a $j_0$, it suffices to notice that by \eqref{27.11.06.2} and b), $j_0=(j'+j-1)/2$ is convenient.

 On Figure \ref{7.12.08.2},  give an illustration of the case $(\gamma_2)$  with $l=4$, $j_0=2$. The comments on this figure are the same as the ones on Figure \ref{7.12.08.1}. 
\begin{figure}\begin{center}$$\begin{array}{ccc}
\xymatrix{&2\ar@/_1pc/@{.}[r]|\Delta\ar[r]|-2&3\ar[dr]|-1\\ 1\ar@{.}[rrr]|\Delta\ar[ur]|-1&&&4\ar[d]|-2\\ 8\ar@{.}[rrr]|\Delta\ar[u]|-2&&&5\ar[dl]|-1\\ &7\ar@/^1pc/@{.}[r]|\Delta\ar[ul]|-1&6\ar[l]|-2}
&&
\xymatrix{
\{2,3\}\ar@(ul,ur)[]|2\ar@/^/[d]^1\\ \{1,4\}\ar@/^/[u]^1\ar@/^/[d]|2\\ \{8,5\}\ar@/^/[d]^1\ar@/^/[u]|2\\ \{7,6\}\ar@(dr,dl)[]|2\ar@/^/[u]^1
} \\ &&
\\ \substack{\textrm{restriction of the graph $G(\sigma, g_1g_2)$ to}\\ \{(l,i, 1), (l,i,2),\ldots, (l,i,8) \} }&& \substack{\textrm{restriction of the graph $G(\sigma, g_1g_2)/\Delta$ to the}\\ \textrm{set of classes of } (l,i, 1), (l,i,2),\ldots, (l,i,8) }
\end{array}$$\caption{}\label{7.12.08.2}\end{center}\end{figure}

 \indent  \indent  ($\gamma_3$) Now, notice that the case ($\gamma_2$) implies the existence of two monochromatic directed cycles with length $1$ in $G(\sigma, g_1g_2)/\Delta$: indeed, the edges $$(l,i,j_0)\to (l,i,c_{2l}(j_0)), \quad  (l,i,c_{2l}^l(j_0))\to (l,i,c_{2l}^{l+1}(j_0))$$ of $G(\sigma, g_1g_2)$ give rise to the edges $$
 \{(l,i,j_0), (l,i,c_{2l}(j_0))\}\xymatrix{\ar@(dr,ur)[] }\quad\quad\quad \{(l,i,c_{2l}^l(j_0)), (l,i,c_{2l}^{l+1}(j_0))\}\xymatrix{\ar@(dr,ur)[] }
 $$ in $G(\sigma, g_1g_2)/\Delta$. The  respective colors of these edges are $$\begin{cases}1,1&\textrm{if $j_0$ is odd and $l$ is  even,}\\ 1,2&\textrm{if $j_0$ is odd and $l$ is odd,}\\  2,2&\textrm{if $j_0$ is even and $l$ is even,}\\ 2,1&\textrm{if $j_0$ is even and $l$ is odd.}\end{cases}$$
 Hence this case is excluded if $A_1=A_2=\{2\}$, and if $A_1=\{2\}$, $A_2=\{1,2\}$, then $j_0$ and $l$ have to be even.

 e) So we have proved that the only non singleton classes of $\Delta$ are of two types
 
 \indent I.  $\{(l,i,j),(l,i',j')\}$
with $l\in [p]$, $i\neq i'\in [N_l(\sigma)]$, $j,j'\in [2l]$, and  if there is such a class   $\{(l,i,j),(l,i',j')\}$ in $\Delta$, then  there is  a unique $q\in [l]$ (depending symmetrically on $i$ and $i'$) \st  the classes, in $\Delta$,   of $(l,i, 1), (l,i,2),\ldots, (l,i,2l) (l,i', 1), (l,i',2),\ldots, (l,i',2l) $ are given by \eqref{26.06.07.1}. We shall denote this integer $q$ by $q_\Delta(l,\{i,i'\})$. 

\indent II. $\{(l,i,j),(l,i,j')\}$ 
with $l\in [p]$, $i\in [N_l(\sigma)]$, $j,j'\in [2l]$, and that if there is such a class   $\{(l,i,j),(l,i,j')\}$ in $\Delta$, then  there is  $j_0\in [2l]$, unique up to a replacement by $c_{2l}^l(j_0)$,  \st  the classes, in $\Delta$,   of $(l,i, 1),  (l,i,2),\ldots, (l,i,2l)$ are given by  (\ref{knives.out.30.11.06}). 
We denote by $T_l(\Delta)$ the set of such $i$'s in $[N_l(\sigma)]$. 
 If $A_1=A_2=\{2\}$, then this case cannot occur, hence $T_l(\Delta)=\emptyset$. 
  If $A_1=A_2=\{1,2\}$, then there is no restriction on $l$ and $j_0$, and  we shall denote the unique element of   $\{j_0, c_{2l}^l(j_0)\}\cap [l]$ by $q_\Delta(l,\{i\})$. 
 If $A_1=\{2\}$, $A_2=\{1,2\}$, then the only restriction on  $l, j_0$ is that they have to be  even, and we shall denote the unique element of   $\{\ff{2}j_0, \ff{2}c_{2l}^l(j_0)\}\cap [\f{l}{2}]$ by $q_\Delta(l,\{i\})$.

 If $i\in [N_l(\sigma)]$ is \st for all $j\in [2l]$, the singleton $\{(l,i,j)\}$ is a class of $\Delta$, then we define $q_\Delta(l,\{i\})$ to be $0$.

f) Let us define, for $l\in [p]$, $P_l(\Delta)$ to be the partition of $[N_l(\sigma)]$ which links to elements $i,i'$ \ssi there is $j,j'\in [2l]$ \st $(l,i,j)$ and $(l,i',j')$ are linked by $\Delta$ (we use the convention that the empty set has a unique partition, that this partition has cardinality zero and that if $N_l(\sigma)=0$,  $P_l(\Delta)$ is this partition).

g) It is easily seen that the number of vertices of $G(\sigma,w)/\Delta$ is $\sum_{l=1}^p2l|P_l(\Delta)|-l|T_l(\Delta)|$, that the sum of the number of $1$-colored   directed cycles of length two  and of  the number of $2$-colored   directed cycles of length two in $G(\sigma,w)/\Delta$ is $\sum_{l=1}^p2l|\{C\in P_l(\Delta)\ste |C|=2\}|+(l-1)|T_l(\Delta)|$ and that the sum of the number of  $1$-colored   directed cycles of length one  and of  the number of $2$-colored    directed cycles of length one in $G(\sigma,w)/\Delta$ is $\sum_{l=1}^p2|T_l(\Delta)|.$  Hence by a), the Neagu characteristic of $G(\sigma,w)/\Delta$ is 
\begin{eqnarray*}&\lf( \sum_{l=1}^p2l|P_l(\Delta)|-l|T_l(\Delta)|\ri)
-2p+\lf(
\sum_{l=1}^p2l|\{C\in P_l(\Delta)\ste |C|=2\}|+(l-1+1)|T_l(\Delta)|
\ri)&\\ 
&= \quad -2p+2\sum_{l=1}^pl(|P_l(\Delta)|+|\{C\in P_l(\Delta)\ste |C|=2\}|).&\end{eqnarray*} But since for all $l$, $P_l(\Delta)$ is a partition of $[N_l(\sigma)]$ where all classes have cardinality one or two, $|P_l(\Delta)|+|\{C\in P_l(\Delta)\ste |C|=2\}|=|[N_l(\sigma)]|=N_l(\sigma)$. Hence since $\sum_{l=1}^plN_l(\sigma)=p$, the Neagu characteristic of $G(\sigma,w)/\Delta$ is 
null.

i) In the case where $A_1=A_2=\{1,2\}$, let us prove that  cardinality of  $C(\sigma, g_1g_2, \{1,2\}, \{1,2\})$ is 
$\prod_{l=1}^p \E[(\sqrt{l}X+l+1)^{N_l(\sigma)}],$ 
for $X$   standard Gaussian variable. Note first that if $C(\sigma)$ denotes $\{l\in [p]\ste N_l(\sigma)\neq 0\}$, then 
 this reduces to  $$\prod_{l\in C(\sigma)} \sum_{i=0}^{ \lfloor N_l(\sigma)/2\rfloor}{N_l(\sigma) \choose 2i}l^i\E(X^{2i})(l+1)^{N_l(\sigma)-2i}.$$ 
Let us
define, for $N\geq 1$, $\Part_{2,1}(N)$ to be the set of partitions of $[N]$ in which all classes have cardinality $1$ or $2$.
It is easy to see that the function \begin{eqnarray*}\vfi \, :\, C(\sigma, g_1g_2, \{1,2\}, \{1,2\})&\to& \{(P_l,(q(l,A))_{A\in P_l})_{l\in C(\sigma)}\ste \forall l\in C(\sigma), P_l\in \Part_{2,1}(N_l(\sigma)),\\ &&\forall A \in P_l\textrm{ \st $|A|=2$, } q(l,A)\in [l] \textrm{ and } \\  && \forall A \in P_l\textrm{ \st $|A|=1$, } q(l,A)\in \{0\}\cup[l]\} \\
\Delta & \mapsto & (P_l(\Delta),(q_\Delta(l,A))_{A\in P_l(\Delta)})_{l\in C(\sigma)}
\end{eqnarray*}is a bijection. Hence 
 it suffices to prove that for all $l\geq 1$, $ N\geq 1$, the cardinality of the set of pairs $(P,(q(A))_{A\in P})$ \st $P\in \Part_{2,1}(N)$ and for all  $A \in P$, $q(A)\in [l]$ if $|A|=2$ and $q(A)\in \{0\}\cup[l]$ in the other case, is equal to $$\sum_{i=0}^{ \lfloor N/2\rfloor}{N\choose 2i}l^i\E(X^{2i})(l+1)^{N-2i}.$$ It follows easily from the well known fact (called Wick's formula) that for all $i\geq 1$, the number of partitions of a set with cardinality $2i$ in which all classes have two elements is equal to $\E(X^{2i})$.

To  prove that  cardinalities of  $C(\sigma, g_1g_2, \{2\}, \{2\})$ and of   $C(\sigma, g_1g_2, \{2\}, \{1,2\})$ are the ones given in the statement of the Lemma, 
we use the same technique, replacing the bijection $\vfi$ by 

\begin{eqnarray*}\phi \, :\, C(\sigma, g_1g_2, \{2\}, \{2\})&\to& \{(P_l,(q(l,A))_{A\in P_l})_{l\in C(\sigma)}\ste \forall l\in C(\sigma), P_l\in \Part_{2,1}(N_l(\sigma))\textrm{ and }\\ &&\forall A \in P_l\textrm{ \st $|A|=2$, } q(l,A)\in [l], \\  && \forall A \in P_l\textrm{ \st $|A|=1$, } q(l,A)=1\}\\
\Delta & \mapsto & (P_l(\Delta),(q_\Delta(l,A))_{A\in P_l(\Delta)})_{l\in  C(\sigma)}
\end{eqnarray*}
\begin{eqnarray*}\psi \, :\, C(\sigma, g_1g_2, \{2\}, \{1,2\})&\to& \{(P_l,(q(l,A))_{A\in P_l})_{l\in C(\sigma)}\ste \forall l\in C(\sigma), P_l\in \Part_{2,1}(N_l(\sigma)),\\ &&\forall l\in C(\sigma), \forall A \in P_l\textrm{ \st $|A|=2$, } q(l,A)\in [l],\\  &&  \forall l\textrm{ even, }\forall A \in P_l\textrm{ \st $|A|=1$, } q(l,A)\in \{0\}\cup[l/2]\})\textrm{ and}\\ &&\forall l\textrm{ odd, }\forall A \in P_l\textrm{ \st $|A|=1$, } q(l,A)=1\}\\
\Delta & \mapsto & (P_l(\Delta),(q_\Delta(l,A))_{A\in P_l(\Delta)})_{l\in C(\sigma)}
\end{eqnarray*}
\end{pr}

\noindent{\bf Proof of Theorem \ref{mardi.28.11.06.2}}. Let us prove (i). (ii) and (iii) can be treated analogously.   
We are going to apply Theorem 1.1 of \cite{fbg.1}.    For $q\geq 1$,  $k_1, \ldots, k_q\geq 1$, for  $p=k_11+\cdots +k_qq$ and $\sigma\in \Sy_{p}$ \st $N_1(\sigma)=k_1$, \ldots, $N_q(\sigma=k_q)$ , by Proposition \ref{propo.les.pates.cuisent.et.beber.attend.10.06} and Lemma  \ref{mardi.28.11.06.1}, $$\f{n^{p}}{1^{k_1}\cdots q^{k_q}k_1!\cdots k_q!}P(\{\forall i=1,\ldots,p, \sigma_n(i)=\sigma(i)\}$$ converges, 
 as $n$ goes to infinity, to $\prod_{l=1}^q \ff{k_l!} \E[((\sqrt{l}X+l+1)/l)^{k_l}]$, 
 where $X$ is a standard Gaussian   variable. Hence  we have to prove that 
 for all $r_1,\ldots, r_q\geq 0$, the series $$\ds\sum_{k_1\geq  r_1}\cdots \sum_{k_q\geq  r_q}(-1)^{k_1-r_1+\cdots +k_q-r_q}{k_1\choose r_1}\cdots {k_q\choose r_q} \prod_{l=1}^q \ff{k_l!} \E[((\sqrt{l}X+l+1)/l)^{k_l}]$$
 converges to $ \prod_{l=1}^q\mu_l(r_l)$, which is equal, by Lemma \ref{1.12.06.apres.tango}, to 
 $ \prod_{l=1}^q e^{-1-\ff{2l}}\f{\E[(X+\sqrt{l})^{r_l}]}{{r_l}!\, l^{r_l/2}}.$ 
Since this series factorize  and  by Remark 1.2 of \cite{fbg.1},    it suffices to prove that  for  $l\geq 1$, $r\geq 1$, the series  
  $\sum_{k\geq r}(-1)^{k-r} {k\choose r} \ff{k!}\E[((\sqrt{l}X+l+1)/l)^{k}]$ 
   converges to $e^{-1-\ff{2l}}\f{\E[(X+\sqrt{l})^r]}{r!\, l^{r/2}}$. 
 This follows from an application of the Dominated Convergence Theorem.\hfill  $\square$

\end{document}